\documentclass[11pt,twoside, final]{amsart}
\copyrightinfo{0}{Iranian Mathematical Society}
\pagespan{1}{\pageref*{LastPage}}
\usepackage{etoolbox,lastpage}
\commby{}
\date{\scriptsize   Received: , Accepted: .}
\usepackage{amsmath,amsthm,amscd,amsfonts,amssymb,enumerate}
\usepackage{graphicx}		
\usepackage{color}
\usepackage[colorlinks]{hyperref}
\usepackage{verbatim} 
\newtheorem{theorem}{Theorem}[section]
\newtheorem{proposition}[theorem]{Proposition}
\newtheorem{lemma}[theorem]{Lemma}
\newtheorem{corollary}[theorem]{Corollary}
\theoremstyle{definition}
\newtheorem{definition}[theorem]{Definition}
\newtheorem{example}[theorem]{Example}
\theoremstyle{remark}

\numberwithin{equation}{section}

\def\a1s{a_1,\cdots, a_s}
\def\a{\alpha}

\def\aa{\mathcal A}

\def\andd{\quad\hbox{and}\quad}

\def\b{\beta}

\def\bb{\mathcal{B}}

\def\bl4{B_{\ell\geq4}}

\def\d{\delta}

\def\End{\hbox{End}}

\def\bbbf{\mathbb{F}}

\def\gg{{\mathcal G}}

\def\fg{\mathfrak{g}}

\def\hh{{\mathcal H}}

\def\fh{\mathfrak{h}}

\def\LL{\mathcal{L}}

\def\ep{\epsilon}
\def\fm{(\cdot,\cdot)}

\def\bbbq{\mathbb{Q}}

\def\1k{\frac{1}{k}}
\def\op{\oplus}
\def\ot{\otimes}

\def\la{\langle}
\def\ra{\rangle}

\def\sub{\subseteq}
\def\sg{\sigma}

\def\rcross{R^{\times}}

\def\pf{\noindent{\bf Proof. }}

\def\rre{R_{re}}
\def\rim{R_{ns}}

\def\fp{\mathfrak{p}}

\def\i{{\mathcal I}}

\def\T{{\mathcal T}}

\def\u{{\mathcal U}}

\def\v{{\mathcal V}}

\def\w{{\mathcal W}}

\def\bbbz{{\mathbb Z}}

\def\1il{1\leq i\leq\ell}

\begin{document}


\title[Locally finite simple  Lie superalgebras]{Locally finite basic classical simple  Lie superalgebras}

\author[M. Yousofzadeh]{Malihe Yousofzadeh$^*$}
\address[Malihe Yousofzadeh]{Department of Mathematics, University of Isfahan, Isfahan, Iran,
P.O.Box 81745-163, and School of Mathematics, Institute for Research in
Fundamental
Sciences (IPM), P.O. Box: 19395-5746, Tehran, Iran}
\email{ma.yousofzadeh@sci.ui.ac.ir $\&$ ma.yousofzadeh@ipm.ir.}

  \thanks{$^*$Corresponding author}
%

 \maketitle
%

\begin{abstract}
In this work, we study direct limits of finite dimensional basic classical simple Lie superalgebras and obtain the conjugacy classes of Cartan subalgebras under the group of automorphisms.\\
\textbf{Keywords:}  Locally finite Lie superalgebra, Finite dimensional basic classical simple Lie superalgebra, Cartan subalgebra, Conjugacy classes.  \\
\textbf{MSC(2010):}  Primary: 17B40; Secondary: 17B65.
\end{abstract}

\section{Introduction}
In 1990, R. H\o egh-Krohn and B. Torresani  \cite{HT} introduced irreducible  quasi simple Lie algebras  as a generalization of both affine Lie algebras and finite dimensional simple Lie algebras over the complex numbers. In 1997, the authors in   \cite{AABGP} systematically studied irreducible  quasi simple Lie algebras under the name extended affine Lie algebras. A nonzero Lie algebra  is called an extended affine Lie algebra if it is equipped with an invariant nondegenerate symmetric bilinear form and that it has a weight space decomposition with respect to a finite dimensional Cartan subalgebra (i.e., a finite dimensional self-centralizing  toral subalgebra)  whose root vectors satisfy some natural conditions. Working with toral subalgebras in place of finite dimensional self-centralizing toral subalgebras, E. Neher \cite{N1} defines the notion of invariant affine reflection algebras. In \cite{you3}, the author introduces and studies the super version of invariant affine reflection algebras called extended affine Lie superalgebras.  Finite dimensional basic classical simple Lie superalgebras and affine Lie superalgebras are examples of extended affine Lie superalgebras having a Cartan subalgebra.

One knows that affine Lie algebras are realized using loop algebras.  In \cite{ABFP}, the authors deal with the realization of extended affine Lie algebras as a generalization of affine Lie algebras. Extended affine Lie algebras has a very close connection with certain  kind of root graded Lie algebras \cite{BM} to which we refer as finite-Lie tori so that the realization of extended affine Lie algebras comes back to the realization of finite-Lie tori. In \cite{ABFP}, the authors using multiloop algebras instead of loop algebras, obtain almost all finite-Lie tori.
The ingredients to construct a multiloop algebra in order to obtain a finite-Lie torus, is a finite dimensional simple Lie algebra $\fg$ and a finite sequence of finite ordered  commuting  automorphisms of $\fg.$

On the other hand, affine Lie superalgebras are obtained using a loop superalgebra starting form a finite dimensional basic classical simple Lie superalgebra \cite{van-de}. Also, it is proved that a simple extended affine Lie superalgebra having a Cartan subalgebra is a direct limit of finite dimensional basic classical simple Lie superalgebras. We call these Lie superalgebras locally finite basic classical simple Lie superalgebras.
In a paper under preparation, we urge to realize extended affine Lie superalgebras; to this end,  we must work with multiloop superalgebras staring from a locally finite basic classical simple Lie superalgebra. For this, we first need to know the structure and the classification of locally finite basic classical simple Lie superalgebras.

In this work, we classify all locally finite basic classical simple Lie superalgebras and then study the conjugacy classes of cartan subalgebras under the group of automorphisms.
Locally finite basic classical simple Lie superalgebras with zero odd part are exactly locally finite split simple Lie algebras which are introduced, studied  and  classified by K-H. Neeb and  N. Stumme in \cite{NS}.

We organize this paper as follows: In Section 1, we gather some preliminaries which  we need throughout the paper.  In Section 2, we introduce some locally finite basic classical simple Lie superalgebras and show that they are mutually non-isomorphic. In Section 3, we classify locally finite basic classical simple Lie superalgebras and study the conjugacy classes of Cartan subalgebras; to do this, we need to know the concept of Chevalley bases for finite dimensional basic classical simple Lie superalgebras. A subsection of Section 3 is exclusively   devoted to Chevalley bases and related topics.

\smallskip

 %

\section{Preliminaries}
Throughout this paper, $\bbbf$ is a field of characteristic zero and  $\bbbz_2:=\{0,1\}$ is the unique abelian group of order $2.$ Unless otherwise mentioned, all vector spaces are considered over $\bbbf.$
We  denote the dual space of a  vector space $V$ by  $V^*.$  We denote the degree of a homogenous element $v$ of a superspace by $|v|$ and make a convention that if  in an expression, we use $|u|$ for an element $u$ of a superspace, by default we have assumed $u$ is homogeneous.    We denote the group of automorphisms of an abelian  group $A$ or a Lie superalgebra $A$   by $Aut(A)$ and use $\simeq$ to show the isomorphism between two algebraic structures. For a subset $S$ of an abelian group, by $\la S\ra,$ we mean the subgroup generated by $S$ and for a set $S,$ by $|S|,$ we mean the cardinal number of $S.$ For a map $f:A\longrightarrow B$ and $C\sub A,$ by $f\mid_{C},$ we mean the restriction of $f$ to $C.$ For two symbols $i,j,$ by $\d_{i,j},$ we mean the Kronecker delta. We  also use  $\biguplus$ to  indicate the disjoint union.
We finally recall that the direct union is, by definition,  the direct limit of a direct system
whose morphisms are inclusion maps.

\smallskip

In the sequel, by a {\it symmetric form} on an additive abelian group $A,$ we mean a map $\fm: A\times A\longrightarrow \bbbf$ satisfying
\begin{itemize}
\item $(a,b)=(b,a)$ for all $a,b\in A,$
\item $(a+b,c)=(a,c)+(b,c)$ and $(a,b+c)=(a,b)+(a,c)$ for all $a,b,c\in A.$
\end{itemize}
In this case, we set  $A^0:=\{a\in A\mid(a,A)=\{0\}\}$ and call it the {\it radical} of the form $\fm.$ The form is called {\it nondegenerate} if $A^0=\{0\}.$
We note that if the form is nondegenerate, $A$ is torsion free and we can identify $A$ as a subset of $\bbbq\ot_\bbbz A.$ Throughout the paper, if an abelian group  $A$ is  equipped with a nondegenerate symmetric form, we consider $A$ as a subset of $\bbbq\ot_\bbbz A$ without further explanation.
Also if $A$ is a vector space over $\bbbf,$ bilinear forms are used in the usual sense.
\medskip

\begin{definition}[{\cite[Def. 1.1]{you6}}]\label{iarr}
{\rm Suppose that $A$ is a nontrivial additive abelian group, $R$ is a subset of $A$ and  $\fm:A\times A\longrightarrow \bbbf$ is  a symmetric form. Set
$$\begin{array}{l}
R^0:=R\cap A^0,\\
\rcross:=R\setminus R^0,\\
\rcross_{re}:=\{\a\in R\mid (\a,\a)\neq0\},\;\;\;\rre:=\rcross_{re}\cup\{0\},\\
\rcross_{ns}:=\{\a\in R\setminus R^0\mid (\a,\a)=0 \},\;\;\; \rim:=\rcross_{ns}\cup\{0\}.
\end{array}$$
We say $(A,\fm,R)$ is an {\it extended affine root supersystem}  if the following hold:
$$\begin{array}{ll}
(S1)& \hbox{$0\in R$ and $\la R\ra= A,$}\\\\
(S2)& \hbox{$R=-R,$}\\\\
(S3)&\hbox{for $\a\in \rre^\times$ and $\b\in R,$ $2(\a,\b)/(\a,\a)\in\bbbz,$}\\\\
(S4)&\parbox{4.5in}{ ({\it root string property}) for $\a\in \rre^\times$ and $\b\in R,$  there are nonnegative  integers  $p,q$  with $2(\b,\a)/(\a,\a)=p-q$ such that \begin{center}$\{\b+k\a\mid k\in\bbbz\}\cap R=\{\b-p\a,\ldots,\b+q\a\};$\end{center}we call  $\{\b-p\a,\ldots,\b+q\a\}$ the {\it $\a$-string through $\b,$}
} \\\\
(S5)&\parbox{4.5in}{for $\a\in \rim$ and $\b\in R$ with $(\a,\b)\neq 0,$
$\{\b-\a,\b+\a\}\cap R\neq \emptyset.$ }
\end{array}$$
If there is no confusion, for the sake of simplicity, we say   {\it $R$ is an extended affine root supersystem in $A.$}
Elements of $R^0$ are called {\it isotropic roots,} elements of  $\rre$ are  called  {\it real roots} and  elements of $\rim$ are  called  {\it nonsingular roots}.
A subset $X$ of $R^\times$ is called {\it connected} if each two elements $\a,\b\in X$ are connected in $X$ in the sense that there is a chain $\a_1,\ldots,\a_n\in X$ with $\a_1=\a,$ $\a_n=\b$ and $(\a_i,\a_{i+1})\neq0,$  $i=1,\ldots,n-1.$ An extended affine root supersystem $R$ is called {\it irreducible} if  $\rre\neq\{0\}$ and $\rcross$ is connected (equivalently, $R^\times$ cannot be written as a disjoint union of two nonempty orthogonal subsets).
An extended affine root supersystem $(A,\fm,R)$ is called a {\it locally finite root supersystem} if the form $\fm$ is nondegenerate and it is called an {\it  affine reflection system} if $\rim=\{0\};$ see \cite{N1}.}
\end{definition}

\begin{definition}
{\rm Suppose that $(A,\fm, R)$ is a locally finite root supersystem.
\begin{itemize}
\item
 The subgroup $\w$ of $Aut(A)$ generated by $r_\a$ ($\a\in\rre^\times$) mapping  $a\in A$ to $a-\frac{2(a,\a)}{(\a,\a)}\a,$  is called the {\it Weyl group} of $R.$
\item A subset $S$ of  $R$ is called a  {\it sub-supersystem} if the restriction of the form to $\la S\ra$ is nondegenerate, $0\in S,$ for $\a\in S\cap\rre^\times, \b\in S$ and $\gamma\in S\cap\rim$ with $(\b,\gamma)\neq 0,$ $r_\a(\b)\in S$ and  $\{\gamma-\b,\gamma+\b\}\cap
    S\neq\emptyset;$   see \cite[Lem. 1.4 \& Rem. 1.6(ii)]{you6}.
\item A sub-supersystem $S$ of $R$ is called  {\it closed} if for $\a,\b\in S$ with $\a+\b\in R,$ we have $\a+\b\in S.$
\item If  $(A,\fm, R)$ is irreducible, $R$ is said to be  of {\it real type} if  $\hbox{span}_\bbbq R_{re}=\bbbq\ot_{\bbbz} A;$ otherwise, we say it is of {\it imaginary type.}
\item The locally finite  root supersystem $(A,\fm,R)$ is called a {\it locally finite root system} if $\rim=\{0\};$ see \cite{LN}.
\item $(A,\fm, R)$ is  said to be {\it isomorphic} to another  locally finite root supersystem $(B,\fm',S)$ if there is a group isomorphism $\varphi:A\longrightarrow B$ and a nonzero scalar $r\in\bbbf$ such that $\varphi(R)=S$ and  $(
a_1,a_2)=r(\varphi(a_1),\varphi(a_2))'$ for all $a_1,a_2\in A.$ In this case, we write $R\simeq S.$
\end{itemize}}
\end{definition}
\begin{lemma}\label{super-sys}
Suppose that $( A,\fm, R)$ is a locally finite root supersystem with Weyl group $ \w.$ Then we have the following:

(i) For $ A_{re}:=\la R_{re}\ra$ and $\fm_{re}:=\fm\mid_{_{ A_{re}\times  A_{re}}},$ $( A_{re},\fm_{re}, R_{re})$ is a locally finite root system.

(ii) If $ R$ is irreducible and $ R_{ns}\neq\{0\},$ then $ R_{ns}^\times=\w\d\cup-\w\d$ for each $\d\in  R_{ns}^\times.$
\end{lemma}
\pf   See \cite[\S 3]{you2}.
\qed

\medskip

Using Lemma \ref{super-sys}, to know the classification of  irreducible locally finite root supersystems, we first need to know the classification of locally finite root systems. Suppose that $T$ is a nonempty  index set  with $|T|\geq 2$ and $\u:=\op_{i\in
T}\bbbz\ep_i$ is the free $\bbbz$-module over   the
set $T.$ Define the  form $$\begin{array}{c}\fm:\u\times\u\longrightarrow\bbbf\\
(\ep_i,\ep_j)\mapsto\d_{i,j}, \hbox{ for } i,j\in T
\end{array}$$
and set
\begin{equation}\label{locally-finite}
\begin{array}{l}
\dot A_T:=\{\ep_i-\ep_j\mid i,j\in T\},\\
D_T:=\dot A_T\cup\{\pm(\ep_i+\ep_j)\mid i,j\in T,\;i\neq j\},\\
B_T:=D_T\cup\{\pm\ep_i\mid i\in T\},\\
C_T:=D_T\cup\{\pm2\ep_i\mid i\in T\},\\
BC_T:=B_T\cup C_T.
\end{array}
\end{equation}
These are irreducible locally finite root systems
in their $\bbbz$-span's. Moreover, each irreducible  locally finite root system is either  an irreducible finite root system  or a locally finite root system  isomorphic to one of these locally finite root
systems. We refer to locally finite root systems listed in (\ref{locally-finite}) as  {\it type} $A,D,B,C$
and $BC$  respectively. We note that if $R$ is  an irreducible locally finite
root system as above, then  $(\a,\a)\in\bbbz^{\geq0}$ for all
$\a\in R.$ This allows us to  define
$$\begin{array}{l}
R_{sh}:=\{\a\in R^\times\mid (\a,\a)\leq(\b,\b);\;\;\hbox{for all $\b\in R$} \},\\
R_{ex}:=R\cap2 R_{sh}\andd
R_{lg}:= R^\times\setminus( R_{sh}\cup R_{ex}).
\end{array}$$
The elements of $R_{sh}$ (resp. $R_{lg},R_{ex}$) are called {\it
short roots} (resp. {\it long roots, extra-long roots}) of $R$. We point  out that following the usual notation in the literature,  the locally finite root system  of type $A$  is denoted by  $\dot A$ instead of $A,$ as all locally finite root systems listed above are spanning sets for $\bbbf\ot_\bbbz \u$ other than the one of type $A$ which spans a subspace of codimension 1; see \cite{LN} and  \cite[Rem. 1.6(i)]{you6}.

\medskip

Next suppose that $I,J$ are two index sets with $I\cup J\neq\emptyset$ and  $F$ is  a free abelian group with a basis $\{\ep_i,\d_j\mid i\in I,j\in J\}.$ Define a form $\fm:F\times F\longrightarrow \bbbf$ with $$(\ep_i,\ep_r):=\d_{i,r}, (\d_j,\d_s):=-\d_{j,s}\andd (\ep_i,\d_j)=0;\;\;\;\; i,r\in I,\;j,s\in J.$$
\newpage
Set
{\footnotesize\begin{equation}\label{types}\begin{array}{rl}
\dot A(I,I):=&\pm\{\ep_i-\ep_r,\d_i-\d_r,\ep_i-\d_r-\frac{1}{\ell}\sum_{k\in I}(\ep_k-\d_k)\mid i,r\in I\};\;\\
 &({\small \ell:=|I|\in\bbbz^{\geq2}}),\\\\
\dot A(I,J):=&\pm\{\ep_i-\ep_r,\d_j-\d_s,\ep_i-\d_j\mid i,r\in I,j,s\in J\};\;\\
 &({\small |I|\neq |J| \hbox{ if $I,J$ are finite sets}}),\\\\
B(I,J):=&\pm\{\ep_i,\d_j,\ep_i\pm\ep_r,\d_j\pm\d_s,\ep_i\pm\d_j\mid i,r\in I,j,s\in J,\;i\neq r\},\\\\
C(I,J):=&\pm\{\ep_i\pm\ep_r,\d_j\pm\d_s,\ep_i\pm\d_j\mid i,r\in I,j,s\in J\},\\\\
D(I,J):=&\pm\{\ep_i\pm\ep_r,\d_j\pm\d_s,\ep_i\pm\d_j\mid i,r\in I,j,s\in J,i\neq r\},\\\\
BC(I,J):=&\pm\{\ep_i,\d_j,\ep_i\pm\ep_r,\d_j\pm\d_s,\ep_i\pm\d_j\mid i,r\in I,j,s\in J\},\\\\
F(4):=&\pm\{0,\ep_1,\d_i\pm\d_j,\d_i,\frac{1}{2}(\ep_1\pm\d_1\pm\d_2\pm\d_3)\mid 1\leq i\neq j\leq 3\};\\
 &({\small I=\{1\},J=\{1,2,3\}}),\\\\
G(3):=&\pm\{0,\ep_1,2\ep_1,\d_i-\d_j,2\d_i-\d_j-\d_k
,\ep_1\pm(\d_i-\d_j)\mid \{i,j,k\}=\{1,2,3\}\};\\
&({\small I=\{1\},J=\{1,2,3\}})
\end{array}\end{equation}}
in which if $I$ or $J$ is empty, the corresponding indices disappear.  We mention that the $\bbbz$-span of all these locally finite  root supersystems are $F$ except for $\dot A(I,J),$ so  to denote this type, we use $\dot A$ instead of $A.$

When $I=\{1\}$ and $J$ is a nonempty index set, we denote  $D(I,J)$ by $C(J),$ so we have
\begin{equation*}
\label{CJ}
C(J)=\{0,\pm\d_j\pm\d_s,\pm\ep_1\pm\d_j\mid j,s\in J\}.
\end{equation*}

In the sequel if either $I$  or  $J$ is a  finite  set, we may replace it by its cardinality in each type, e.g., we may denote $B(I,J)$ by $B(|I|,|J|)$ if $I$ and $J$ are finite sets. Using this convention,  $C(1)$ can be identified with $\dot A(1,2)$ and   $D(2,1)$ is nothing but the root system of the  finite dimensional basic classical simple Lie superalgebra  $D(2,1,\a)$ for $\a=1.$
We drew the attention of readers to the point that  our notations  have a  minor difference with the notations in the literature, more precisely,  $C(n)$ for $n\in\bbbz^{\geq 1}$ and $\dot A(m,n)$ for $m,n\in\bbbz^{\geq 1}$ in our sense are denoted by  $C(n+1)$ and $A(m-1,n-1)$ respectively in the literature. Our  notations allow us  to  switch smoothly from the finite case  to the infinite case.

\begin{theorem}[{\cite[\S 4]{you2}}]\label{classification}
{ Each irreducible locally finite root supersystem is  either isomorphic to  the root system of  a finite dimensional basic classical simple Lie superalgebra or  isomorphic to one of the  root supersystems introduced in (\ref{types}). Among all irreducible locally finite root supersystems,  $C(I)$ and $\dot A(I,J)$} with $|I|\neq |J|$ if both $I$ and $J$ are finite, are  of imaginary type and the other ones are of real type.
\end{theorem}

\begin{lemma}\label{lin}
Suppose that $\v$ is a vector space equipped with a symmetric bilinear form $\fm$ and $R$ is a subset of $\v$ such that $(\la R\ra ,\fm_{\la R\ra\times \la R\ra}, R)$ is a locally finite root supersystem. Suppose that  $\{\a_1,\ldots,\a_n\}\sub \rre$ is $\bbbq$-linearly independent. Then $\{\a_1,\ldots,\a_n\}$ is $\bbbf$-linearly independent; also if  $\rim\setminus \hbox{span}_\bbbq \{\a_1,\ldots,\a_n\}\neq \emptyset$ and $\d\in \rim\setminus \hbox{span}_\bbbq \{\a_1,\ldots,\a_n\},$ then $\{\d,\a_1,\ldots,\a_n\}$ is also $\bbbf$-linearly independent.
\end{lemma}
\begin{proof} We assume
$\rim\setminus \hbox{span}_\bbbq \{\a_1,\ldots,\a_n\}\neq \emptyset$ and $\d\in \rim\setminus \hbox{span}_\bbbq \{\a_1,\ldots,\a_n\},$ and show  $\{\d,\a_1,\ldots,\a_n\}$ is  $\bbbf$-linearly independent; the other statement is similarly proved. Take  $\{1,x_i\mid i\in I\}$  to be  a basis for $\bbbq$-vector space $\bbbf.$
Suppose that $r,r_1,\ldots,r_n\in\bbbf$ and $r\d+\sum_{j=1}^nr_j\a_j=0.$  Suppose that  for $1\leq j\leq n,$ $r_j=s_j+\sum_{i\in I}s_j^ix_i$ with $\{s_j,s_j^i\mid i\in I\}\sub \bbbq.$ We first show  $r=0.$ To the contrary, assume $r\neq 0.$ Without loss of generality, we assume $r=1.$ So $0=\d+\sum_{j=1}^n r_j\a_j=\d+\sum_{j=1}^n(s_j+\sum_{i\in I}s_j^ix_i)\a_j.$ Now for $\a\in \rre,$ we have
$$\frac{2(\d,\a)}{(\a,\a)}+\sum_{j=1}^ns_j\frac{2(\a_j,\a)}{(\a,\a)}+\sum_{j=1}^n\sum_{i\in I}s_j^ix_i\frac{2(\a_j,\a)}{(\a,\a)}=0.$$ This implies that  for $\a\in\rre$ and  $i\in I,$ $\sum_{j=1}^ns_j^i\frac{2(\a_j,\a)}{(\a,\a)}=0$ and so $(\sum_{j=1}^ns_j^i\a_j,\a)=0.$ But it follows from Lemma \ref{super-sys}$(b)(i)$ that  the form on $\hbox{span}_\bbbq \rre$ is nondegenerate, so $\sum_{j=1}^ns_j^i\a_j=0$ for all $i\in I.$ Now as $\{\a_j\mid 1\leq j\leq n\}$  is $\bbbq$-linearly independent, we have
 $$
s_j^i=0 \;\;\;\;(i\in I,j\in \{1,\ldots,n\}).
$$
Therefore, we get $0=\d+\sum_{j=1}^n r_j\a_j=\d+\sum_{j=1}^n(s_j+\sum_{i\in I}s_j^ix_i)\a_j=\d+\sum_{j=1}^ns_j\a_j.$ Thus we have $\d=-\sum_{j=1}^ns_j\a_j$ which is absurd. This shows that $r=0.$ Now repeating the above argument, one gets that $s_j^i=0$  for all $ i\in I,j\in \{1,\ldots,n\}$  and that $0=\sum_{j=1}^ns_j\a_j.$ Thus we have $s_j=0$ for all $1\leq j\leq n.$ This  implies that $r_j=s_j+\sum_{i\in I}s_j^ix_i=0$ for all $1\leq j\leq n$ and so we are done.
\end{proof}

\begin{lemma}[{\cite[Lem. 2.3]{you6}}]
Suppose that $R$ is an irreducible locally finite root supersystem of type $X$ in an  abelian group $A.$ Then we have the following:

(i) $A$ is a  free abelian group and $R$ contains a   $\bbbz$-basis for $A.$

(ii) If $X\neq A(\ell,\ell),$ $R$ contains a $\bbbz$-basis  $\Pi$  for $A$ satisfying  the partial sum property in the sense that     for each $\a\in \rcross,$ there are $\a_{1},\ldots,\a_{n}\in\Pi$  (not necessarily distinct) and $r_1,\ldots,r_n\in\{\pm1\}$ with $\a=r_1\a_{1}+\cdots+r_n\a_{n}$ and  $r_1\a_{1}+\cdots+r_t\a_{t}\in \rcross,$  for all $1\leq t\leq n.$
 \end{lemma}

\begin{definition}
{\rm A subset $\Pi$  of a locally finite root supersystem  $R$ is called an {\it integral base} for $R$ if $\Pi$ is a $\bbbz$-basis for $A.$ An integral base $\Pi$ of $R$ is called a {\it base} for $R$ if it satisfies the partial sum property.}
\end{definition}

\begin{lemma}[{\cite[Lem. 2.4($iii$)]{you6}}]\label{base}
If $R$ is an infinite irreducible locally finite root supersystem, then there is a base $\Pi$ for $R$ and a class $\{R_\gamma\mid \gamma\in\Gamma\}$ of finite irreducible  closed sub-supersystems of $R$ of the same type as $R$ such that $R$ is the direct union of $R_\gamma$'s and for each $\gamma\in\Gamma,$ $\Pi\cap R_\gamma$ is a base for $R_\gamma.$
\end{lemma}

\section{Locally finite basic classical simple Lie superalgebras}\label{Locally finite basic classical simple Lie superalgebras}
We recall that a Lie superalgebra $\gg$ is called {\it locally finite} if each finite subset of $\gg$ generates a finite dimensional subsuperalgebra.
Suppose that $\LL=\LL_0\op\LL_1$ is a nonzero Lie superalgebra equipped with a nondegenerate invariant even supersymmetric bilinear form $(\cdot,\cdot)$ and $\hh$ is a nontrivial subalgebra  of $\LL_{0}$ such that with respect to $\hh,$ $\LL$ has a weight space decomposition $\LL=\op_{\a\in\hh^*}\LL^\a$ via the adjoint representation and the restriction of the form  $\fm$ to $\hh$ is nondegenerate.
We call $R:=\lbrace \alpha\in \mathfrak{h}^{\ast} \mid \LL^{\alpha}\neq \{0 \}\rbrace$, the \textit{root system} of $\LL$ (with respect to $\mathfrak{h}$). Each element of  $R$ is called a \textit{root}. We mention that    $\mathfrak{h}$ is abelian and as  $\LL_{0}$ as well as $\LL_{1}$ are $\mathfrak{h}$-submodules of $\LL$, we have using \cite[Pro. 2.1.1]{MP} that $\LL_{0}=\bigoplus_{\alpha\in \mathfrak{h}^{\ast}} \LL_{0}^{\alpha}$ and $\LL_{1}=\bigoplus_{\alpha\in \mathfrak{h}^{\ast}} \LL_{1}^{\alpha}$ with $\LL_{i}^{\alpha}:=\LL_{i}\cap \LL^{\alpha},  i=0, 1.$ We refer to the elements of $R_{0}:=\lbrace \alpha\in \mathfrak{h}^{\ast} \mid  \LL_{0}^{\alpha}\neq \{0\} \rbrace$ (resp. $R_{1}:=\lbrace \alpha \in \mathfrak{h}^{\ast} \mid  \LL_{1}^{\alpha}\neq \{0\}\rbrace$) as \textit{even roots} (resp. \textit{odd roots}) and  note that $R=R_{0}\cup R_{1}$.
Since the form is invariant and even,  for $\a,\b\in R$ and $ i,j\in\{0,1\},$  we have  \begin{equation*}\label{form-zero}(\LL_i^\a,\LL_{j}^\b)=\{0\}\;\;\; \;\; \hbox{if $i\neq j$ or $\a+\b\neq0$}.\end{equation*} This   in particular implies that for $i\in\{0,1\}$ and $\a\in R_i,$   the form restricted to  $\LL_i^\a+\LL_i^{-\a}$ is nondegenerate.
 Take  $\mathfrak{p}:\hh\longrightarrow \hh^*$ to be  the function  mapping   $h\in\hh$ to $(h,\cdot).$    Since the form is nondegenerate on $\hh,$ the map $\mathfrak{p}$ is one to one (and so onto if $\hh$ is finite dimensional). So for each element $\a$ of the image $\hh^\frak{p}$ of $\hh$ under the map $\mathfrak{p},$  there is a unique $t_\a\in\hh$ representing $\a$ through the form $\fm.$ Now we can transfer the form on $\hh$ to a form on $\hh^\mathfrak{p},$ denoted again by $\fm$  and defined by \begin{equation*}\label{form}(\a,\b):=(t_\a,t_\b)\;\;\;(\a,\b\in \hh^\fp).\end{equation*}
Using  Lemma 3.1 of \cite{you3}, if $\a\in R\cap \hh^\mathfrak{p},$  $x\in\LL^\a$ and $y\in\LL^{-\a}$ with  $[x,y]\in\hh,$    we have \begin{equation*}\label{bracket}[x,y]=(x,y)t_\a.\end{equation*}
  We also  draw the attention of readers to the point that if either $\hh$ is finite dimensional or $\LL^0=\hh,$ it is not hard to see that $R\sub \hh^{\frak{p}}.$

\begin{definition}
{\rm A  Lie superalgebra $\LL=\LL_{0}\op\LL_{1},$   is called a {\it locally finite basic classical simple Lie superalgebra} if
\begin{itemize}
\item $\LL$ is   locally finite and  simple,
\item  $\LL$ is equipped with an invariant nondegenerate even supersymmetric  bilinear form.
\item $\LL_{0}$ has a nontrivial subalgebra $\hh$ (refereed to as a {\it Cartan subalgebra}) with respect to which $\LL$ has a weight space decomposition $\LL=\sum_{\a\in \hh^*}\LL^\a$ via the adjoint representation with corresponding root system $R$ such that $\LL^0=\hh$ and $\rcross=\{\a\in R\mid(\a,R)\neq\{0\}\}\neq \emptyset.$
    \end{itemize} We may also write $(\LL,\hh,\fm)$ is a locally finite basic classical simple Lie superalgebra.
}
\end{definition}

\begin{theorem}[{\cite[Thm. 2.30]{long-ver}}]\label{main1}
Suppose that $(\LL,\hh,\fm)$ is a locally finite basic classical simple  Lie superalgebra, then

(i) the root system $R$ of $\LL$ is an irreducible locally finite root supersystem,

(ii) $\LL$ is a direct union of finite dimensional  basic classical simple  Lie  superalgebras,

(iii) $[\LL_{0},\LL_{0}]$ is a semisimple Lie algebra,

(iv) if $\LL_{1}\neq\{0\},$ it  is a completely reducible $\LL_{0}$-module with at most two irreducible constituents.
\end{theorem}

In the rest of this section, we shall introduce some non-isomorphic examples of locally finite basic classical simple Lie superalgebras. Let us start with some notations.
For a unital associative superalgebra $\aa$ and nonempty index supersets $I,J,$ by an $I\times J$-matrix
with entries in $\aa,$ we mean a map $A:I\times J\longrightarrow \aa.$ For $i\in I,j\in J,$ we set
$a_{ij}:=A(i,j)$ and call it the {\it $(i,j)$-th entry} of $A.$ By a convention, we denote the matrix
$A$ by $(a_{ij}).$  We also denote the set of all $I\times J$-matrices  with entries in $\aa$ by
$\aa^{I\times J}.$  If $I=J,$ we denote $\aa^{I\times J}$ by $\aa^I.$  For $A=(a_{ij})\in\aa^{I\times J},$ the matrix
$B=(b_{ij})\in\aa^{J\times I}$ with  $$b_{ij}:=\left\{ \begin{array}{ll} a_{ji}& |i|=|j|\\
a_{ji}& |i|=1,|j|=0\\
-a_{ji}& |i|=0,|j|=1
\end{array}
\right.$$  is called the
{\it supertransposition} of $A$ and denoted by $A^{st}.$
 If $A=(a_{ij})\in\aa^{I\times J}$ and
$B=(b_{ij})\in \aa^{J\times K}$ are such that for all $i\in I$ and $k\in K,$ at most for finitely many $j\in J,$ $a_{ij}b_{jk}$'s are nonzero, we define the product $AB$ of $A$ and $B$ to be the
$I\times K$-matrix $C=(c_{ik})$ with $c_{ik}:=\sum_{j\in J}a_{ij}b_{jk}$ for all $i\in I,k\in K.$
We note that if $A,B,C$ are three matrices such that $AB,$ $(AB)C,$ $BC$ and $A(BC)$ are defined,
then $A(BC)=(AB)C.$ We make a convention that if $I$ is a disjoint union of subsets
$I_1,\ldots, I_t$ of $I,$ then for an $I\times I$-matrix $A,$ we write
$$A=\left [\begin{array}{ccc} A_{1,1}&\cdots& A_{1,t}\\
A_{2,1}&\cdots&A_{2,t}\\
\vdots&\vdots&\vdots\\
A_{t,1}&\cdots&A_{t,t}\\
\end{array}\right ]$$
in which for $1\leq r,s\leq t,$ $A_{r,s}$ is an $I_r\times I_s$-matrix whose $(i,j)$-th entry
coincides with $(i,j)$-th entry of $A$ for all $i\in I_r,j\in I_s.$ In this case, we say that $A\in \aa^{I_1\uplus\cdots\uplus I_t}$ and note that the defined matrix product  obeys the product of block matrices. If $\{a_i\mid i\in I\}\sub \aa,$ by
$\hbox{diag}(a_i),$ we mean an $I\times I$-matrix whose $(i,i)$-th entry is $a_i$ for all $i\in I$ and other
entries are zero. If $\aa$ is unital, we set $1_I:=\hbox{diag}(1_\aa).$ A matrix $A\in\aa^I$ is called {\it invertible} if there is a matrix $B\in\aa^I$ such that $AB$ as well as  $BA$ are defined and $AB=BA=1_I;$ such a $B$ is unique and denoted by $A^{-1}.$
 For $i\in I,j\in J$ and $a\in \aa,$ we define $E_{ij}(a)$ to be a matrix in $\aa^{I\times J}$ whose $(i,j)$-th entry is
$a$ and other entries are zero and if $\aa$ is unital, we set $$e_{i,j}:=E_{i,j}(1).$$  Take $M_{I\times J}(\aa)$ to be the subspace of $\aa^{I\times J}$ spanned by
$\{E_{ij}(a)\mid i\in I,j\in J,a\in A\}.$ $M_{I\times J}(\aa)$  is a superspace with $M_{I\times J}(\aa)_i:=\hbox{span}_\bbbf\{E_{r,s}(a)\mid |r|+|s|+|a|=\bar i\},$ for $i=0,1.$ Also with respect to the  multiplication of matrices, the vector superspace
$M_{I\times I}(\aa)$ is an associative $\bbbf$-superalgebra and so it is a Lie superalgebra under the Lie bracket
$[A,B]:=AB-(-1)^{|A||B|}BA$ for all $A,B\in M_{I\times I}(\aa).$ We denote  this Lie superalgebra  by
$\mathfrak{pl}_I(\mathcal{A}).$ For $X,Y\in \mathfrak{pl}_I(\mathcal{A}),$ we have  $(XY)^{st}=(-1)^{|X||Y|}Y^{st}X^{st}.$
 Finally, for an element $X\in \mathfrak{pl}_I(\mathcal{A}),$ we set $str(X):=\sum_{i\in I}(-1)^{|i|}x_{i,i}$ and call it the {\it supertrace} of $X.$

\begin{lemma}\label{iso1}
(i) Suppose that $Q$ is a homogeneous element of $\bbbf^I,$ then $\gg_Q:=\{X\in \mathfrak{pl}_{I}(\bbbf)\mid X^{st}Q=-(-1)^{|X||Q|}QX\}$ is a Lie subsuperalgebra of $\mathfrak{pl}_{I}(\bbbf).$

(ii) If $Q_1,Q_2$ are homogeneous elements of $\bbbf^I$  and $T$ is an invertible homogeneous element of $\bbbf^I$ of degree zero such that   $Q_2=T^{st}Q_1 T,$ then $\gg_{Q_1}$ is isomorphic to $\gg_{Q_2}$ via the isomorphism mapping $X$ to $T^{-1}XT.$

(iii) Suppose that $I$ and $J$ are two supersets and $\eta:I\longrightarrow J$ is a bijection preserving the degree. For a matrix $A=(A_{ij})$ of $\bbbf^I,$ define $A^\eta\in\bbbf^J$ to be $(A^\eta_{ij})$ with $A^\eta_{ij}=A_{\eta^{-1}(i)\eta^{-1}(j)}.$ If $Q$ is a homogeneous element of $\bbbf^I$ and $Q':=Q^\eta,$ then the Lie superalgebra $\gg_Q:=\{X\in \mathfrak{pl}_I(\bbbf) \mid X^{st}Q=-(-1)^{|X||Q|}QX\}$ is isomorphic to the Lie superalgebra  $\gg_{Q'}:=\{X\in \mathfrak{pl}_J(\bbbf)\mid X^{st}Q'=-(-1)^{|X||Q'|}Q'X\}.$
\end{lemma}

\begin{proof}
$(i),(ii)$ It is easy to verify.

$(iii)$ Suppose that  matrices $A,B\in\bbbf^I$ are  such that $AB$ is defined, then for $i,j\in I,$  we have
\begin{eqnarray*}
(A^\eta B^\eta)_{\eta(i)\eta(j)}=\sum_{t\in I}A^\eta_{\eta(i)\eta(t)}B^\eta_{\eta(t)\eta(j)}=\sum_{t\in I}A_{it} B_{tj}=(AB)_{ij}=(AB)^\eta_{\eta(i)\eta(j)}.
 \end{eqnarray*}
This in particular  implies that if $A,B,C,D\in\bbbf^I$ are such that $AB$ and $CD$ are defined and $AB=CD,$ then $A^\eta B^\eta=C^\eta D^\eta.$
 Moreover, as $\eta$ preserves the degree, we have $(A^{st})^\eta=(A^\eta)^{st}.$ Now it is easy to see that the function $\theta:\gg_Q\longrightarrow \gg_{Q'}$ mapping $X$ to $X^\eta$ is an isomorphism.
\end{proof}

\begin{example}\label{exa1}
{\rm
For two disjoint  index sets  $I,J$ with $J\neq\emptyset,$ suppose that $\{0,i,\bar i,\mid i\in I\cup J\}$ is a superset
with $|0|=|i|=|\bar i|=0$ for $i\in I$ and $|j|=|\bar j|=1$ for $j\in J.$  We set $\dot I:=I\cup\bar I,$ $\dot I_0:=\{0\}\cup I\cup\bar I$ and $\dot J:=J\cup \bar J$ where $$\bar I:=\{\bar i\mid i\in I\}\andd \bar J:=\{\bar j\mid j\in J\}.$$ For $\mathcal{I}=\dot I\cup\dot J$ or $\mathcal{I}=\dot I_0\cup\dot J,$ we set $$Q_\mathcal{I}:=\left(\begin{array}{cc}M_1&0\\
0&M_2\end{array}\right)$$ in which $$ M_2:=\sum_{j\in J}(e_{j,\bar j}-e_{\bar j, j})\;\&\;M_1:=\left\{\begin{array}{ll}
-2e_{0,0}+\sum_{i\in I}(e_{i,\bar i}+e_{\bar i,i})& \hspace{-2mm}\hbox{if }  \hbox{\small $\mathcal{I}=\dot I_0\cup\dot J$}\\
\sum_{i\in I}(e_{i,\bar i}+e_{\bar i, i})& \hspace{-2mm}\hbox{if } \hbox{\small $I\neq\emptyset, \mathcal{I}=\dot I\cup\dot J.$}
\end{array}\right.$$
Now by Lemma \ref{iso1}, $$\gg_\i:=\gg_{Q_\i}=\{X\in \mathfrak{pl}_\i(\bbbf)\mid X^{st}Q_\i=-Q_\i X\}$$ is a Lie subsuperalgebra of $\mathfrak{pl}_\i(\bbbf)$ which  we refer to as  $\mathfrak{osp}(2I,2J)$ or $\mathfrak{osp}(2I+1,2J)$ if  $\i=\dot I\cup\dot J$ or  $\i=\dot I_0\cup\dot J$ respectively.
Set\begin{equation*}\label{h}\fh:=\hbox{span}_\bbbf\{h_t,d_k\mid t\in I,\;k\in J\}\end{equation*} in which for  $t\in I$ and $ k\in J,$ $$h_t:=e_{t,t}-e_{\bar t,\bar t}\andd d_k:=e_{k,k}-e_{\bar k,\bar k}$$ and for $i\in I$ and $j\in J,$ define
$$\begin{array}{lll}
\begin{array}{c}
\ep_i:\fh\longrightarrow \bbbf\\
h_t\mapsto \d_{i,t},\;\;\; d_k\mapsto 0,
\end{array}&&
\begin{array}{c}
\d_j:\fh\longrightarrow \bbbf\\
h_t\mapsto 0,\;\;\; d_k\mapsto \d_{j,k},
\end{array}
\end{array}$$  $(t\in I,k\in J).$ One sees that   $\gg_\i$ has a weight space decomposition with respect to $\fh.$
Taking $R(\i)$  to be the corresponding  set of weights,  we have
{\small \begin{equation*}\label{root}\begin{array}{l}R(\dot I_0\cup\dot J)=\{\pm\ep_r,\pm(\ep_r\pm\ep_s),\pm\d_p,\pm(\d_p\pm\d_q),\pm(\ep_r\pm\d_p)\mid  r, s\in I,p,q\in J,r\neq s\},\\
\\R(\dot I\cup\dot J)=\{\pm(\ep_r\pm\ep_s),\pm(\d_p\pm\d_q),\pm(\ep_r\pm\d_p)\mid  r, s\in I,\;p,q\in J,\;r\neq s\}\end{array}\end{equation*}} in which $\pm(\ep_r\pm\ep_s)$'s
 are disappeared if $|I|=1$ and $\pm\ep_r$'s, $\pm(\ep_r\pm\ep_s)$'s as well as $\pm(\ep_r\pm\d_p)$'s are disappeared if $|I|=0.$
 Moreover, for $r, s\in I,\;p,q\in J$  with $r\neq s$ and $p\neq q,$ we have
{\small$$\begin{array}{ll}
(\gg_\i)^{\ep_r+\ep_s}=\hbox{span}_\bbbf  (e_{r,\bar s}- e_{s,\bar r}),&
(\gg_\i)^{-\ep_r-\ep_s}=\hbox{span}_\bbbf  (e_{\bar r,s}- e_{\bar s,r}),\\\\
(\gg_\i)^{\ep_r-\ep_s}=\hbox{span}_\bbbf  (e_{r,s}- e_{\bar s,\bar r)},&
(\gg_\i)^{\d_p+\d_q}=\hbox{span}_\bbbf (e_{p,\bar q}+ e_{q,\bar p}),\\\\
(\gg_\i)^{-\d_p-\d_q}=\hbox{span}_\bbbf (e_{\bar p,q}+ e_{\bar q,p}),&
(\gg_\i)^{\d_p-\d_q}=\hbox{span}_\bbbf (e_{p, q}- e_{\bar q,\bar p}),\\\\
(\gg_\i)^{\ep_r+\d_p}=\hbox{span}_\bbbf (e_{r,\bar p}+ e_{p,\bar r}),&
(\gg_\i)^{-\ep_r-\d_p}=\hbox{span}_\bbbf (e_{\bar r,p}- e_{\bar p,r}),\\\\
(\gg_\i)^{\ep_r-\d_p}=\hbox{span}_\bbbf (e_{r,p}- e_{\bar p,\bar r}),&
(\gg_\i)^{-\ep_r+\d_p}=\hbox{span}_\bbbf (e_{\bar r,\bar p}+ e_{p,r}),\\\\
(\gg_\i)^{2\d_p}=\hbox{span}_\bbbf e_{p,\bar p},&
(\gg_\i)^{-2\d_p}=\hbox{span}_\bbbf e_{\bar p,p},\\\\
\mathfrak{osp}(2I+1,2J)^{\ep_r}=\hbox{span}_\bbbf  ( e_{0,\bar r}+2e_{r,0}),&
\mathfrak{osp}(2I+1,2J)^{-\ep_r}=\hbox{span}_\bbbf  (e_{0,r}+2e_{\bar r,0}),\\\\
\mathfrak{osp}(2I+1,2J)^{\d_p}=\hbox{span}_\bbbf ( e_{0,\bar p}-2e_{ p,0}),&
\mathfrak{osp}(2I+1,2J)^{-\d_p}=\hbox{span}_\bbbf (e_{0,p}+2 e_{\bar p,0}).
\end{array}$$}
Define $$\fm:\gg_{\i}\times\gg_\i\longrightarrow \bbbf;\;\;\; (x,y)\mapsto str(xy)\;\;\;\;\;(x,y\in\gg_\i).$$
Then $(\gg_\i,\fh,\fm)$ is a locally finite basic classical simple Lie superalgebra whose root system is an irreducible locally finite root supersystem  of type $X$ as in the following table:

\begin{center}
\begin{tabular}{|c|c|c||c|c|c|}
\hline
$X$&$(|I|,|J|)$&$\i$&$X$&$(|I|,|J|)$&$\i$\\
\hline
$B(0,J)$& $(0,\geq1)$&$\dot I_0\cup\dot J$&$ C(J)$ & $(1,\geq2)$& $\dot I\cup\dot J$\\
\hline
$B(1,J)$&$(1,\geq 1)$&$\dot I_0\cup\dot J$&$D(2,1,\a)$ & $(2,1)$& $\dot I\cup\dot J$\\
\hline
$B(I,1)$ & $(\geq2,1)$&$\dot I_0\cup\dot J$&$D(2,J)$ & $(2,\geq2)$& $\dot I\cup\dot J$\\
\hline
$B(I,J)$&$(\geq 2,\geq2)$&$\dot I_0\cup\dot J$&$D(1,I)$ & $(\geq3,1)$&$\dot I\cup\dot J$\\
\hline
$\dot A(0,2)$&$(1,1)$&$\dot I\cup\dot J$&
$D(I,J)$ & $(\geq3,\geq2)$&$\dot I\cup\dot J$\\
\hline
\end{tabular}
\end{center}
We refer to $\fh$ as the {\it standard Cartan subalgebra} of $\gg_\i.$ We note that $(\gg_\i)_{0}$ is centerless unless  $\i=\dot I\cup \dot J$ with  $| I|=1;$ see Lemma 2.33 of \cite{long-ver}. In this case, suppose $I=\{1\},$ then  for a  fixed index  $j\in J,$ $t_{\ep_1+\d_j}-(1/2)t_{2\d_j}$ is a nonzero central element of the even part of $\gg_\i.$
\hfill{$\diamondsuit$}
}
\end{example}

As in \cite[\S 1]{NS}, we have the following lemma:
\begin{lemma}
\label{cor2}
Suppose that $I,J$ are two nonempty  index sets with $|I|=\infty,$ then $\mathfrak{osp}(2I,2J)\simeq\mathfrak{osp}(2I+1,2J).$
\end{lemma}

\begin{proof}
 Consider the following matrices of $\bbbf^{\{0\}\uplus I\uplus \bar I\uplus J\uplus\bar J}:$
{\scriptsize$$S:=\left (\begin{array}{ccccc}
1&0&0&0&0\\
0&I&I&0&0\\
0&I&-I&0&0\\
0&0&0&I&0\\
0&0&0&0&I
\end{array}\right), Q_e:=\left (\begin{array}{ccccc}
-2&0&0&0&0\\
0&2I&0&0&0\\
0&0&-2I&0&0\\
0&0&0&0&I\\
0&0&0&-I&0
\end{array}\right),Q:=\left (\begin{array}{ccccc}
-2&0&0&0&0\\
0&0&I&0&0\\
0&I&0&0&0\\
0&0&0&0&I\\
0&0&0&-I&0
\end{array}\right),$$}
then we  have $S^{st} Q S=Q_e.$ Also for matrices  {\scriptsize$$S':=\left (\begin{array}{cccc}
I&I&0&0\\
I&-I&0&0\\
0&0&I&0\\
0&0&0&I
\end{array}\right), Q_o:=\left (\begin{array}{cccc}
2I&0&0&0\\
0&-2I&0&0\\
0&0&0&I\\
0&0&-I&0
\end{array}\right),Q':=\left (\begin{array}{cccc}
0&I&0&0\\
I&0&0&0\\
0&0&0&I\\
0&0&-I&0
\end{array}\right)$$} of $\bbbf^{I\uplus \bar I\uplus J\uplus\bar J},$
we have   ${S'}^{st} Q' S'=Q_o.$ Now by Lemma \ref{iso1}, $\gg_{Q}\simeq\gg_{Q_e},$  $\gg_{Q'}\simeq\gg_{Q_o},$ and   $\gg_{Q_e}\simeq\gg_{Q_o}.$ This completes the proof.
\end{proof}

\begin{example}\label{exa2}
{\rm
Suppose that $J$ is a superset  with $J_{0},J_{1}\neq\emptyset$ \textcolor[rgb]{1.00,0.00,0.00}{and $(|J_0|,|J_1|)\neq(1,1)$}. Set $\gg:=\mathfrak{sl}(J_{0},J_{1})=\{X\in \mathfrak{pl}_J(\bbbf)\mid str(X)=0\}$ and  $\hh:=\hbox{span}_\bbbf\{e_{i,i}-e_{r,r},e_{j,j}-e_{s,s},e_{i,i}+e_{j,j}\mid i,r\in J_{0},j,s\in J_{1} \}.$
For $t\in J_{0}, k\in J_{1},$ define $$\begin{array}{ll}\ep_t:\hh\longrightarrow \bbbf,\\
\hbox{\small$e_{i,i}-e_{r,r}\mapsto \d_{i,t}-\d_{r,t},\;e_{j,j}-e_{s,s}\mapsto 0,e_{i,i}+e_{j,j}\mapsto \d_{i,t},$}\\
\\\d_k:\hh\longrightarrow \bbbf,\\
\hbox{\small$e_{i,i}-e_{r,r}\mapsto 0,\;e_{j,j}-e_{s,s}\mapsto \d_{j,k}-\d_{k,s},e_{i,i}+e_{j,j}\mapsto \d_{k,j},$}\end{array}
$$ ($i,r \in J_{0}, j,s\in J_{1}$). Also define $$\fm:\gg\times\gg\longrightarrow \bbbf;\;\; (X,Y)\mapsto str(XY).$$ If $|J|<\infty$ and $|J_{0}|=|J_{1}|,$ take  $K:=\bbbf \sum_{j\in J}e_{jj}$  and note that it is a subset of the radical of the form $\fm.$ So it induces a bilinear form on $\gg/K$ denoted again by $\fm.$  Set $$\mathfrak{sl}_s(J_{0},J_{1}):=\left\{\begin{array}{ll}
\gg/K&\hbox{if $|J|<\infty$ and $|J_{0}|=|J_{1}|$}\\
\gg& \hbox{otherwise}.
\end{array}\right.$$
Then $(\LL:=\mathfrak{sl}_s(J_{0},J_{1}),\fm,\hh/K)$ is a locally finite basic classical simple Lie superalgebra with root system $$\{\ep_i-\ep_j,\d_p-\d_q,\pm(\ep_i-\d_p)\mid i,j\in J_{0}, p,q\in J_{1}\}$$  which is an irreducible  locally finite root supersystem of type $X$ as in the following table:

\begin{center}
\begin{tabular}{|c|c|}
\hline
$X$&$(|J_{0}|,|J_{1}|)$\\
\hline
$\dot A(0,J_{1})$& $(1,\geq2)$\\
\hline
$\dot A(0,J_{0})$& $(\geq2,1)$\\
\hline
$\dot A(J_{0}, J_{1})$&
$(\geq2,\geq2)$\\
&
$|J_{0}|\neq |J_{1}| \hbox{ if $J_{0},J_{1}$ are both finite }$\\
\hline
$A(\ell,\ell)$ & $(\ell,\ell)$\;\;\; $(\ell\in\bbbz^{\geq1})$\\
\hline
\end{tabular}
\end{center}
Also if \textcolor[rgb]{1.00,0.00,0.00}{ $(|J_{0}|,|J_{1}|)\neq (2,2),$} for $ i,j\in J_{0}$ and $ p,q\in J_{1}$ with $i\neq j$ and $p\neq q,$ we have
$$
\begin{array}{ll}
\LL^{\ep_i-\ep_j}=\bbbf e_{i,j},& \LL^{\d_p-\d_q}=\bbbf e_{q,p},\\
\LL^{\ep_i-\d_p}=\bbbf e_{i,p},& \LL^{-\ep_i+\d_p}=\bbbf e_{p,i}.
\end{array}
$$
We refer to $\hh/K$ as the {\it standard Cartan subalgebra} of $\LL=\mathfrak{sl}_s(J_{0},J_{1}).$ We now  need to discuss the center of $\LL_{0}$ for our future purpose. We recall from finite dimensional theory of Lie superalgebras that  if $|J_{0}|, |J_{1}|<\infty,$ $\LL_{0}$ has nontrivial  center if and only if $|J_{0}|\neq |J_{1}|$ and that in this case, it has a one dimensional center.  Now suppose  $|J_{0}\cup J_{1}|=\infty,$ say $|J_{0}|=\infty.$ Fix $i_0\in J_{0}$ and $j_0\in J_{1}.$ Then $\{e_{i,i}-e_{i_0,i_0},e_{j,j}-e_{j_0,j_0},e_{i_0,i_0}+e_{j_0,j_0}\mid i\in J_{0}\setminus\{i_0\},j\in J_{1}\setminus\{j_0\}\}$ is  a basis for $\hh.$ Suppose $i_1,\ldots,i_\ell$ are distinct elements of $J_{0}\setminus\{i_0\}$ and $j_1,\ldots,j_n$ are distinct elements of
$J_{1}\setminus\{j_0\}.$ If $z=\sum_{t=1}^\ell r_t(e_{i_t,i_t}-e_{i_0,i_0})+\sum_{t=1}^n s_t(e_{j_t,j_t}-e_{j_0,j_0})+k(e_{i_0,i_0}+e_{j_0,j_0})$  (where $\sum_{t=1}^n s_t(e_{j_t,j_t}-e_{j_0,j_0})$ is disappeared if $|J_{1}|=1$)
 is an element of the center of $\LL_{0},$ then for each $i\in J_{0}\setminus\{i_0\},$  $[z,(\LL_{0})^{\ep_{i}-\ep_{i_0}}]=\{0\}.$ Now if $i=i_s$ for some  $s\in \{1,\ldots,\ell\},$ we get   $r_s+(\sum_{t=1}^\ell r_t)-k=0$ and if $i\not\in \{i_1,\ldots,i_\ell\},$ we get $\sum_{t=1}^\ell r_t-k=0.$ Therefore we have $r_s=0$ for all $s\in\{1,\ldots,\ell\}$ and so $k=0.$ This shows that $\LL$ is centerless  if $|J_{1}|=1.$ If  $|J_{1}|>1,$ $[z,\LL^{\ep_{i_0}-\d_{j_0}}]=\{0\}.$ This implies  that $\sum_{t=1}^n s_t=0 .$ We also have  $[z,\LL^{\d_{j}-\d_{j_0}}]=\{0\}$ for all $j\in J_{1}\setminus\{j_0\}.$ Now it follows that $s_t=0$ for all $t\in \{1,\ldots,n\}.$ This means that $z=0$ and so $\LL$ is centerless.
\hfill{$\diamondsuit$}
}
\end{example}

\begin{lemma}\label{class1}
For  index sets  $I,J$ with $|J|\neq 0$ and a superset $T$ with $|T_{0}|,|T_{1}|\neq0,$  set $$\mathfrak{a}_{I,J}:=\mathfrak{osp}(2I,2J) (\hbox{if $I\neq\emptyset$}), \mathfrak{b}_{I,J}:=\mathfrak{osp}(2I+1,2J),\mathfrak{c}_T:=\mathfrak{sl}_s(T_{0},T_{1}).$$
 Suppose that $\gg$ and $\LL$ are two Lie superalgebras of the class $$\{\mathfrak{a}_{I,J},\mathfrak{b}_{I,J},\mathfrak{c}_T\mid I,J,T\}.$$ Then $\gg$ and $\LL$ are isomorphic if and only if (up to changing the role of $\gg$ and $\LL$) one of the following holds:
 \begin{itemize}
\item there are index sets $I,J,I',J'$ with $|I|=|I'|\neq0,$ $|J|=|J'|\neq0,$ $\gg=\mathfrak{a}_{I,J}$ and $\LL=\mathfrak{a}_{I',J'},$
\item there are index sets $I,J,I',J'$ with   $|I|=|I'|,$ $|J|=|J'|\neq0,$ $\gg=\mathfrak{b}_{I,J}$ and $\LL=\mathfrak{b}_{I',J'},$
\item there are supersets $I,J$ with $|I_{0}|=|J_{0}|\neq0,$ $|I_{1}|=|J_{1}|\neq0,$ or $|I_{0}|=|J_{1}|\neq0,$ $|I_{1}|=|J_{0}|\neq0$ such that $\gg=\mathfrak{c}_{I}$ and $\LL=\mathfrak{c}_{J},$
\item there are   index sets $I,J$ with $|I|=|J|=1$ and a superset $T$ with $|T_{0}|=1, |T_{1}|=2$  or $|T_{0}|=2, |T_{1}|=1$ such that $\gg=\mathfrak{a}_{I,J}$ and $\LL=\mathfrak{c}_T,$
\item there are index sets $I,J$ with $J\neq\emptyset,|I|=\infty,$ $\gg=\mathfrak{a}_{I,J}$ and $\LL=\mathfrak{b}_{I,J}.$
 \end{itemize}
Moreover, in each of  the first three cases, the mentioned isomorphism can be chosen such that the standard Cartan subalgebra of $\gg$ is mapped  to the standard Cartan subalgebra of $\LL.$
\end{lemma}
\begin{proof}
We first note that for two Lie algebras  $\mathfrak{k}_1$ and $\mathfrak{k}_2$  such that  $[\mathfrak{k}_1,\mathfrak{k}_1]$ and $[\mathfrak{k}_2,\mathfrak{k}_2]$ are semisimple with the complete  sets  of simple ideals $\{\mathfrak{k}_1^1,\ldots,\mathfrak{k}_1^n\}$ and $\{\mathfrak{k}_2^1,\ldots,\mathfrak{k}_2^m\}$ respectively,  if $\mathfrak{k}_1$ and $\mathfrak{k}_2$ are isomorphic, we have
\begin{equation}\label{fact}
\begin{array}{ll}
\bullet& \hbox{\small$[\mathfrak{k}_1,\mathfrak{k}_1]$ and $[\mathfrak{k}_2,\mathfrak{k}_2]$ are isomorphic,}\\
\bullet&\hbox{\small$\mathfrak{k}_1$ is centerless if and only if $\mathfrak{k}_2$ is centerless,}\\
\bullet& \hbox{\small$m=n$ and (under a permutation of indices) $\mathfrak{k}_1^i\simeq\mathfrak{k}_2^i$  for $i\in\{1,\ldots, n\}.$}
\end{array}
\end{equation}Now take $\aa$ to be one of the Lie superalgebras $\mathfrak{a}_{I,J},\mathfrak{b}_{I,J},\mathfrak{c}_T.$ We have already seen that if  $\aa$  is  infinite dimensional, then the  even part of $\aa$ is centerless if and only if $\aa\neq\mathfrak{a}_{I,J}$ for some infinite index set $J$ and an index set $I$ with $|I|=1.$
Next  suppose that $\gg$ and $\LL$ are as in the statement and assume they are  isomorphic, then  we have $\gg_{0}\simeq \LL_{0}.$ We also know that  $[\gg_{0},\gg_{0}]$ as well as $[\LL_{0},\LL_{0}]$ are semisimple Lie algebras by Theorem \ref{main1}. Using these together with (\ref{fact}),  Lemmas \ref{iso1} and \ref{cor2}, classification of basic classical simple Lie superalgebras  and \cite[Pro.'s VI4,VI6]{NS}, we are done.
\end{proof}

\section{Classification Theorem}
In this section, we classify locally finite basic classical simple Lie superalgebras  (l.f.b.c.s Lie superalgebras for short) and study the conjugacy classes of their Cartan subalgebras under the group of automorphisms.
The first step towards the classification of l.f.b.c.s Lie superalgebras is finding out an isomorphism theorem. One knows that l.f.b.c.s Lie superalgebras with zero odd part are exactly locally finite split simple Lie algebras in the sense of \cite{NS} and that finite dimensional basic classical simple Lie superalgebras and consequently finite dimensional simple Lie algebras are examples of l.f.b.c.s. Lie superalgebras. We know  form  the finite dimensional theory of Lie algebras that due to the interaction of a finite dimensional simple  Lie algebra with its root system, the theorem stating that finite dimensional simple Lie algebras with isomorphic root systems, are isomorphic \cite[Thm. 14.2]{Hum}, plays a crucial role to get the classification of finite dimensional simple Lie algebras. Using this theorem together with the fact that locally finite split simple Lie algebras are a direct union of finite dimensional simple subalgebras, the authors in \cite{NS} prove that two locally finite split simple Lie algebras with isomorphic root systems are isomorphic.  Moreover, they introduce two isomorphic locally finite split simple Lie algebras with non-isomorphic  Cartan subalgebras and isomorphic  root systems.  They use this to find the conjugacy classes of Cartan subalgebras of locally finite split simple Lie algebras.

To get the  classification of l.f.b.c.s Lie superalgebras,  we also prove that  two l.f.b.c.s. Lie superalgebras with isomorphic root systems are isomorphic.
To this end, we  first need to prove the theorem for finite dimensional case. Because of the existence of self-orthogonal roots for a finite dimensional basic classical simple Lie superalgebra, the proof of the mentioned theorem  in the  super case is different from the one in non-super case; more precisely,  we first need to define Chevalley bases for finite dimensional basic classical simple Lie superalgebras.
Chevalley bases for finite dimensional basic classical simple Lie superalgebras  were introduced in 2011 by K. Iohara and  Y. Koga  \cite{KY} using  the fact that a  finite dimensional basic classical simple Lie superalgebra is a contragredient  Lie superalgebra and its Cartan matrix is symmetrizable.  Our definition of Chevalley bases are  somehow different from the one  defined in \cite{KY}.

The zero part of a locally finite basic classical simple Lie superalgebra which is infinite dimensional and not a Lie algebra is either a locally finite split simple Lie algebra or  a direct sum of two locally finite split simple Lie algebras. In the last theorem of this section, we use the result of \cite{NS} to find the conjugacy classes of Cartan subalgebras of locally finite basic classical simple Lie superalgebras.

\begin{lemma}\label{lin2}
Suppose that $(\gg_1,\fm_1,\hh_1),(\gg_2,\fm_2,\hh_2)$ are two locally finite  basic classical simple Lie superalgebras with  corresponding root systems $R_1,R_2$ respectively. For $i=1,2,$ denote the induced form on $\hbox{span}_\bbbf R_i\sub \hh_i^*$ again by $\fm_i.$ If  $R_1$ and $R_2$ are isomorphic, say via $f:\la R_1\ra\longrightarrow \la R_2\ra$ with $( f(\a), f(\a'))_2=k(\a,\a')_1$ for all $\a,\a'\in R$ and some $k\in \bbbf\setminus\{0\},$ then there is a linear isomorphism $\tilde f:\hbox{span}_\bbbf R_1\longrightarrow \hbox{span}_\bbbf R_2$ whose restriction to $\la R_1\ra$ coincides with  $f$ and for $\a,\a'\in \hbox{span}_\bbbf R_1,$ $(\tilde f(\a),\tilde f(\a'))_2=k(\a,\a')_1.$
\end{lemma}
\begin{proof}
We know that $R_1$ is of real type if and only if for each nonsingular root $\d,$ there exists a nonzero integer $n$ with $n\d\in \la (R_1)_{re}\ra$ or equivalently $\hbox{span}_\bbbq (R_1)_{re}=\hbox{span}_\bbbq R_1.$ Now  fix a basis $\{\a_i\mid i\in I\}\sub (R_1)_{re}$ for $\hbox{span}_\bbbq (R_1)_{re}$ as well as a nonzero nonsingular root  $\d$ of $R_1$ if $R_1$ is of imaginary type. Set $$B:=\left\{\begin{array}{ll}
\{\a_i\mid i\in I\} & \hbox{if $R_1$ is of real type}\\
\{\d,\a_i\mid i\in I\} & \hbox{if $R_1 $ is of imaginary type.}
\end{array}\right.$$ Then by Lemma \ref{lin}, $B$ is $\bbbf$-linearly independent and so by Lemma \ref{super-sys}$(b)(ii)$, it is a basis for both $\hbox{span}_\bbbf R_1$ and $\hbox{span}_\bbbq R_1.$ Similarly, $f(B)$ is a basis for $\hbox{span}_\bbbf R_2.$ We define the linear transformation  $\tilde f$ mapping $\a\in B$ to $f(\a).$ It is immediate that $(\tilde f(\a),\tilde f(\a'))=k(\a,\a')$ for $\a,\a'\in \hbox{span}_\bbbf R_1.$ Now if $\a\in R_1\sub\hbox{span}_\bbbq B,$  $\a=\sum_{j=1}^n\frac{r_j}{s_j}\b_j$ where $r_1,s_1,\ldots,r_n,s_n\in\bbbz$ and $\b_1,\ldots,\b_n\in B,$ so for $s=s_1\cdots s_n$ and $r'_j=r_js/s_j$ ($1\leq j\leq n$), we have $s\a=\sum_{j=1}^nr'_j\b_j.$ Therefore, we have $sf(\a)=\sum_{j=1}^nr'_jf(\b_j).$ Thus, we have $f(\a)=\tilde{f}(\a).$
\end{proof}
\subsection{Chevalley bases  for basic classical simple Lie superalgebras}
Suppose that  $\gg$ is a finite dimensional basic classical simple Lie superalgebra of type $X\neq A(1,1)$  with a Cartan subalgebra $\hh$ and  corresponding root system $R=R_0\cup R_1$ such that $\gg_{1}\neq\{0\}.$
In what follows for   $\a\in R^\times$ with $\gg_i^\a\neq \{0\}$ ($i\in\{0,1\}$),  we set  $|\a|:= i.$ Now we want to define a total ordering on  $\v:=\hbox{span}_\bbbq R.$ We fix a basis $\{v_1,\ldots,v_m\}$ for $\v.$ For $u=r_1v_1+\cdots+r_mv_m\in\v,$ we say $0\prec u$ if $u\neq0$ and that the first nonzero $r_i,$ $1\leq i\leq m,$ is positive; next for $u,v\in\v,$ we say $u\prec v$ if $0\prec v-u.$ We set $R^+:=R\cap \{v\in\v\mid 0\prec v\}$ as well as $R^-:=-R^+.$ Elements of $R^+$ are called {\it positive} and elements of $R^-$ are called {\it negative}. As usual,  for $u,v\in\v,$ we say $u\preceq v$ if either $u=v$ or $u\prec v.$ Fix  an invariant  nondegenerate   even supersymmetric  bilinear  form $\fm$ on $\gg.$ We denote the induced  nondegenerate symmetric bilinear   form on $\hh^*$ again by $\fm.$   We recall that  for $\a\in\hh^*,$ $t_\a$ indicates  the unique element of $\hh$ representing $\a$ through the form $\fm.$ For  $\a\in \hh^*,$ set
 $$\sg_\a:=\left\{\begin{array}{ll}
  -1 & \a\in R_1\cap R^-\\
  1& \hbox{otherwise.}
\end{array}\right.$$
Next fix $r\in \bbbf\setminus\{0\}$  and for each $\a\in R^\times,$  set
$$ h_\a:=rt_\a.$$
One can see that
\begin{equation*}
\label{identity}
\sg_{-\a}=(-1)^{|\a|}\sg_\a\andd h_\a=-h_{-\a}\;\;\;(\a\in R^\times).
\end{equation*}
Fixing $Y_{\a}\in \gg^{\a}$ and $Y_{-\a}\in\gg^{-\a}$  with  $[Y_\a,Y_{-\a}]=h_\a$ for $\a\in R^+,$  we have $[Y_\a,Y_{-\a}]=\sg_\a h_\a$ $(\a\in R^\times).$
\begin{definition}
{\rm
A set $\{X_\a,h_{i}\mid \a\in \rcross,i=1,\ldots,\ell\}$ is called a {\it Chevalley basis} for $\gg$ if
\begin{itemize}
\item there are a nonzero scalar $r$ and a   subset $\{\b_1,\ldots,\b_\ell\}$ of $R^\times$ such that $\{h_1:=h_{\b_1},\ldots,h_\ell:=h_{\b_\ell}\}$ is a basis for $\hh$ where for $\a\in R^\times,$ by $h_\a,$ we mean $rt_\a,$
\item for each $\a\in R^\times,$ $X_{\a}\in \gg^{\a},$
\item for each $\a\in R^\times,$ $[X_{\a},X_{-\a}]=\sg_\a h_\a.$
\end{itemize}
}\end{definition}

Suppose that $\{X_\a,h_{i}\mid \a\in \rcross,i=1,\ldots,\ell\}$ is a {\it Chevalley basis} for $\gg.$ We know from Lemma 2.4 of \cite{long-ver} that  if $\a,\b\in \rcross$ such that $\a+\b\in \rcross,$ then $[\gg^\a,\gg^\b]\neq \{0\}.$ This together with the fact that $\dim(\gg^{\a+\b})=1$ implies that there is a nonzero scalar $N_{\a,\b}$ with $[X_\a,X_\b]=N_{\a,\b}X_{\a+\b};$ we also interpret $N_{\a,\b}$ as zero for $\a,\b\in R^\times$ with  $\a+\b\not\in R.$ We refer to $\{N_{\a,\b}\mid \a,\b\in R^\times\}$ as a set of {\it structure constants} for $\gg$ with respect to $\{X_\a,h_{i}\mid \a\in \rcross,i=1,\ldots,\ell\}$.
\begin{proposition}\label{constant}
Keep the same notation as above;  we have the following:

(i) If $\a,\b\in \rcross,$ then  $N_{\a,\b}=-(-1)^{|\a||\b|}N_{\b,\a}.$

(ii)  If $\a,\b\in \rcross$ with $\a+\b\in \rcross,$ then for $s_{\a,\b}:=\sg_\a\sg_{\a+\b},$ we have   $$N_{\a,\b}=s_{\a,\b}N_{\b,-\a-\b}=\sg_\a\sg_{\a+\b}N_{\b,-\a-\b}.$$

(iii) Suppose that $\a,\b\in R^\times$ with  $\a+\b\in R^\times,$ then  $$N_{\a,\b}N_{-\a,-\b}=r_{\a,\b}:=\sg_{\b}\sg_{\b+\a}\sg_\a(-1)^{|\b||\a|}\sum_{i=0}^p(-1)^{i|\a|}(\b-i\a)(h_\a);$$ where $p=0$ if $\a,\b\in \rim$ and  otherwise, $p$ is the largest nonnegative integer such that  $\b-p\a\in R.$

(iv) If $\a,\b,\gamma,\d\in R^\times$ with $\a+\b+\gamma+\d=0$ such that each pair is not the opposite  of the one  another, then {\small$$(-1)^{|\a||\gamma|}\sg_{\a+\b} N_{\a,\b}N_{\gamma,\d}+(-1)^{|\a||\b|}\sg_{\b+\gamma} N_{\b,\gamma}N_{\a,\d}+(-1)^{|\b||\gamma|}\sg_{\a+\gamma} N_{\gamma,\a}N_{\b,\d}=0.$$}
\end{proposition}
\begin{proof}
Using a modified argument as in \cite[Pro. 7.1]{ca}.
\end{proof}

\medskip
We  know  that there are  roots $\a,\gamma$ such that $\a\neq\pm\gamma$ and  $(\a,\gamma)\neq0.$ So either $\a+\gamma\in \rcross$ or $\a-\gamma\in \rcross.$ Replacing $\gamma$ with $-\gamma$ if necessary, we assume $\eta:=-(\a+\gamma)\in\rcross.$ Since $\a+\gamma+\eta=0,$   either two of $\a,\gamma,\eta$ are positive or two of $-\a,-\gamma,-\eta$ are positive.
Selecting this pair of positive roots in an appropriate order, we get a  pair $(\eta_1,\eta_2)$ among the 12 pairs
$$\begin{array}{ll}(\a,\gamma),(\a,\eta),(\gamma,\eta),(\gamma,\a),(\eta,\a),(\eta,\gamma),\\
(-\a,-\gamma),(-\a,-\eta),(-\gamma,-\eta),(-\gamma,-\a),(-\eta,-\a),(-\eta,-\gamma)
\end{array}$$
such that $0\prec \eta_1\preceq \eta_2;$  following \cite{ca}, we call such a pair a {\it special pair}. More precisely, a pair $(\a,\b)$ of elements of $R^\times$ is called a {\it special pair} if $0\prec\a\preceq\b$ and $\a+\b\in R.$  A special pair $(\a,\b)$  is called an {\it extraspecial pair} if for each  special pair $(\d,\gamma)$ with $\a+\b=\d+\gamma,$ we get  $\a\preceq\d.$
\begin{lemma}\label{arbitrary}
Suppose that $\aa$ is the set of all extraspecial pairs  $(\a,\b)$ of $R^\times$  and $\{N_{\a,\b}\mid (\a,\b)\in \aa\}$ is an arbitrary set of nonzero scalars. Then there is $\{e_\a\in\gg^\a\setminus\{0\}\mid \a \in R^+\}$ such that $[e_\a,e_\b]=N_{\a,\b}e_{\a+\b}$ for all $(\a,\b)\in \aa.$
\end{lemma}
\begin{proof}
Suppose that $R^+=\{\a_1,\ldots,\a_n\}$ with $\a_1\prec\ldots\prec\a_n$ and take  $t$ to be   the smallest index such that $\a_t$ is the summation of the components of an extraspecial pair. We choose  arbitrary elements  $e_{\a_i}\in \gg^{\a_i},$ for $1\leq i\leq t-1.$ We know that there is a unique extraspecial pair $(\a,\b)$ with $\a_t=\a+\b,$ so there is a unique pair $(i,j)$ with $i\leq j<t$ such that $\a_t=\a_i+\a_j$ and define $e_{\a_t}=N_{\a_i,\a_j}^{-1}[e_{\a_i},e_{\a_j}].$ Now using an induction process, we can complete the proof; indeed, suppose that $t<r\leq n$  and that $\{e_{\a_s}\mid 1\leq s\leq r-1\}$  with the desired property has been chosen. If $\a_r$ is not the summation of the components of an extraspecial pair, we choose $e_{\a_r}$ arbitrary, but otherwise we pick the unique pair $(i',j')$ with $i'\leq j'<r-1$ such that $\a_r=\a_{i'}+\a_{j'}.$ Now we define  $e_{\a_r}=N_{\a_{i'},\a_{j'}}^{-1}[e_{\a_{i'}},e_{\a_{j'}}].$ This completes the proof.
\end{proof}

\begin{theorem}\label{com-cheval}
Suppose that $\gg$ and $\LL$ are two finite dimensional  basic classical simple Lie superalgebras with Cartan subalgebras $\hh$ and $T$ and corresponding root systems $R=R_0\cup R_1$ and $S=S_0\cup S_1$ respectively which are not of type $A(1,1).$  Suppose that $\fm$ $($resp. $\fm')$ is an invariant  nondegenerate   even supersymmetric bilinear form on $\gg$ $($resp. $\LL)$ and denote the induced forms on $\hh^*$ and $T^*$ again by $\fm$ and $\fm'$ respectively. Suppose that $(\la R\ra ,\fm,R)$ and $(\la S\ra,\fm',S)$ are isomorphic finite root supersystems, say via $f:\la R\ra\longrightarrow \la S\ra.$ Then we have the following:

(i) There are Chevalley bases $\{h_i, e_\a\mid\a\in \rcross,1\leq i\leq \ell\}$ and $\{t_i,x_{\b}\mid \b\in S^\times, 1\leq i\leq \ell\}$ for $\gg$ and $\LL$ with corresponding sets of structure constants $\{N_{\a,\b}\mid \a,\b\in \rcross\}$ and $\{M_{\gamma,\eta}\mid \gamma,\eta\in S^\times\}$ respectively  such that $N_{\a,\b}=M_{f(\a),f(\b)}$ for all $\a,\b\in \rcross.$

(ii) $\{N_{\a,\b}\mid \a,\b\in \rcross\}$ is completely  determined in terms of $N_{\a,\b}$'s for extraspecial pairs $(\a,\b).$

(iii) There is an isomorphism from $\gg$ to $\LL$ mapping $\hh$ to $\T$ and  $e_\a$ to $x_{f(\a)}$ for all $\a\in R\setminus\{0\}.$
\end{theorem}
\begin{proof}
$(i),(ii)$ Suppose that  $k\in\bbbf\setminus\{0\}$ is such that $(f(\a),f(\b))'=k(\a,\b)$ for $\a,\b\in R.$ Fix $r,s\in\bbbf\setminus\{0\}$ such that $r=sk.$ This implies that $r(\a,\b)=sk(\a,\b)=s(f(\a),f(\b))'$ for all $\a,\b\in R.$ Use Lemma \ref{lin2} to extend the map $f$ to a linear isomorphism, denoted again by $f,$ from $\hh^*=\hbox{span}_\bbbf R$ to $T^*=\hbox{span}_\bbbf S$ with
\begin{equation}
\label{compatible}
r(\a,\b)=sk(\a,\b)=s(f(\a),f(\b))'\;\;\;\;\; (\a,\b\in \hh^*).
\end{equation}
For $\a\in \hh^*,$ take $t_\a$ to be the unique element of $\hh$ representing $\a$ through $\fm$ and for $\b\in T^*,$ take $t'_\b$ to be the unique element of $T$ representing $\b$ through $\fm'.$ Next set
\begin{equation*}
\label{h alpha}
h_\a:=rt_\a\andd h'_\b:=st'_\b\;\;\;\;(\a\in R,\;\; \b\in S).
\end{equation*}
Fix a total  ordering $``\preceq"$ on $\hbox{span}_\bbbq R$ as at the beginning of  this subsection and transfer it through $f$ to a total  ordering, denoted again by  $``\preceq",$ on $\hbox{span}_\bbbq S.$ For $\a\in \hh^*$ and $\b\in T^*,$ set
\begin{equation*}
\label{sigma}
\sg_\a:=\left\{\begin{array}{ll}
-1& \hbox{if $\a\in R^-\cap R_1$}\\
1& \hbox{otherwise}
\end{array}\right.\andd \sg'_\b:=\left\{\begin{array}{ll}
-1& \hbox{if $\b\in S^-\cap S_1$}\\
1& \hbox{otherwise.}
\end{array}\right.
\end{equation*}
Suppose that $\aa$ is the set of all extraspecial pairs of $R,$ then
{\small$$\{(f(\a),f(\b))\mid (\a,\b)\in \aa\}=\{(\eta,\gamma)\mid (\eta,\gamma)\hbox{ is an extraspecial pair of $S$}  \}.$$}  Fix  a subset $\{N_{\a,\b}\mid (\a,\b)\in\aa\}$ of nonzero scalars and set $M_{f(\a),f(\b)}:=N_{\a,\b},$ for all $(\a,\b)\in\aa.$ Using Lemma \ref{arbitrary}, one can find $\{e_\a\in\gg^\a\setminus\{0\}\mid\a\in R^+\}$ and $\{x_{\b}\in\LL^\b\setminus\{0\}\mid \b\in S^+\}$ such that $$[e_\a,e_\b]= N_{\a,\b}e_{\a+\b}\andd [x_{f(\a)},x_{f(\b)}]=M_{f(\a),f(\b)}x_{f(\a)+f(\b)};\;\;\;(\a,\b)\in\aa.$$
Now for each $\a\in R^+$ and $\gamma\in S^+,$ choose $e_{-\a}\in\gg^{-\a}$ and $x_{-\gamma}\in \LL^{-\gamma}$ such that $$[e_\a,e_{-\a}]= h_\a\andd [x_\gamma,x_{-\gamma}]= h'_\gamma\;\;\;\;\; (\a\in R^+,\gamma\in S^+)$$ and note that we have
$$
[e_\a,e_{-\a}]=\sg_\a h_\a\andd [x_\gamma,x_{-\gamma}]=\sg'_\gamma h'_\gamma\;\;\;\;\; (\a\in R^\times,\gamma\in S^\times).
$$

Now for each pair $(\a,\b)$ of $\rcross$ with $\a+\b\in \rcross$ and $(\a,\b)\not \in \aa,$ take $N_{\a,\b}$ to be the   unique  nonzero element of $\bbbf$ with $[e_\a,e_\b]=N_{\a,\b}e_{\a+\b};$ also for each pair $(\gamma,\eta)$ of $S^\times$ with $\gamma+\eta\in S^\times$ such that $(\gamma,\eta)$ is not an extraspecial pair, take  $M_{\gamma,\eta}$ to be  the unique nonzero element of $\bbbf$ with $[x_\gamma,x_\eta]=M_{\gamma,\eta}e_{\gamma+\eta}.$
Fix $\{\b_1,\ldots,\b_\ell\}$ such that $\{h_i:=h_{\b_i}\mid 1\leq i\leq \ell\}$ is a basis for $\hh$ and set $t_i:=h'_{f(\b_i)}.$ Then $\{h_i, e_\a\mid\a\in \rcross,1\leq i\leq \ell\}$ and $\{t_i,x_{\b}\mid \b\in S^\times, 1\leq i\leq \ell\}$ are Chevalley bases  for $\gg$ and $\LL$ respectively.
Now contemplating Proposition \ref{constant} and using the same argument as in \cite[Pro. 7.4]{ca}, we get the result.

$(iii)$ Use the same notation as above. Define $\theta:\gg\longrightarrow \LL$ mapping $h_i=h_{\b_i}$ to $t_i=h'_{f(\b_i)}$ and $e_\a$ to $x_{f(\a)}$ for all $\a\in\rcross$ and $1\leq i\leq \ell.$ We claim that $\theta$ is a Lie superalgebra isomorphism. We first note that by  \cite[Pro. 2.5]{long-ver} and \cite[Pro. 3.10]{you3}, $f(R_0)=S_0$ and $f(R_1)=S_1.$  Therefore, we have $\theta(\gg_i)\sub \LL_i$ for $i=0,1.$ Now we need to show $\theta[x,y]=[\theta(x), \theta(y)]$ for all $x,y\in \gg.$ If $x=h_{\b_i}$ and $y=e_\a,$ for some $1\leq i\leq \ell$ and $\a\in \rcross,$ by (\ref{compatible}), we have
\begin{eqnarray*}\theta[h_{\b_i},e_\a]=\theta(\a(h_{\b_i})e_\a)=\a(h_{\b_i})\theta(e_\a)&=&f(\a)(h'_{f(\b_i)})x_{f(\a)}\\
&=&[h'_{f(\b_i)},x_{f(\a)}]\\
&=&[\theta(h_{\b_i}),\theta(e_\a)].\end{eqnarray*} Next suppose  $\a,\b\in \rcross.$ If $\a+\b\not\in R,$ then $f(\a)+f(\b)\not\in S$ and so   $[e_\a,e_\b]=0$ and $[\theta(e_\a),\theta(e_\b)]=[x_{f(\a)},x_{f(\b)}]=0,$ also if $\a+\b\in\rcross,$ then by part ($i$),
\begin{eqnarray*}\theta[e_\a,e_\b]=\theta(N_{\a,\b}e_{\a+\b})=N_{\a,\b}\theta(e_{\a+\b})&=&M_{f(\a),f(\b)}x_{f(\a+\b)}\\
&=&[x_{f(\a)},x_{f(\b)}]\\
&=&[\theta(e_\a),\theta(e_\b)].\end{eqnarray*}
Finally, for $\a\in\rcross,$ if $h_\a=\sum_{i=1}^\ell r_ih_{\b_i}$ for some $r_i\in \bbbf$ $(1\leq i\leq \ell),$ we get $\a=\sum_{i=1}^\ell r_i\b_i$ and so $f(\a)=\sum_{i=1}^\ell r_if(\b_i)$ which in turn implies that \begin{equation}\label{final5}h'_{f(\a)}=\sum_{i=1}^\ell r_ih'_{f(\b_i)}=\sum_{i=1}^\ell r_i\theta(h_{\b_i})=\theta(h_\a).\end{equation} Therefore, we have
\begin{eqnarray*}\theta[e_\a,e_{-\a}]=\theta(\sg_\a h_\a)=\sg_\a h'_{f(\a)}=\sg_{f(\a)} h'_{f(\a)}&=&[x_{f(\a)},x_{-f(\a)}]\\
&=&[\theta(e_\a),\theta(e_{-\a})].
\end{eqnarray*}
This completes the proof.
\end{proof}

Suppose that $\gg$ is a finite dimensional basic classical simple Lie superalgebra with a Cartan subalgebra $\hh$ and corresponding root system $R.$ For a group homomorphism $\phi:\la R\ra\longrightarrow \bbbf\setminus\{0\},$ the linear transformation    $\tilde\phi:\gg\longrightarrow \gg$ mapping $x\in\gg^\a$ ($\a\in R$) to $\phi(\a)x$ is a superalgebra automorphism.

\begin{lemma}\label{cor}
Keep the same notations and assumptions as in Theorem \ref{com-cheval} and its proof. Suppose that $\Pi$ is an integral base for $R$ and fix nonzero elements  $f_\a\in\gg^\a$ and $y_{\a}\in \LL^{f(\a)}$ for all $\a\in \Pi.$ Then there is  an isomorphism  from $\gg$ to $\LL$ mapping $f_\a$ to $y_\a$ and $h_\a$ to $h'_{f(\a)}$ for all $\a\in\Pi.$  Moreover, if $\Pi$ is a base, then such an  isomorphism  is unique.
\end{lemma}
\begin{proof}
Consider the Chevalley bases $\{e_\a,h_i\mid \a\in R^\times,1\leq i\leq \ell \}$ and $\{x_\b,t_i\mid \b\in S^\times,1\leq i\leq \ell\}$  as well as the isomorphism $\theta:\gg\longrightarrow \LL$ as in Theorem \ref{com-cheval}. Since $f_\a\in \gg^{\a}=\bbbf e_\a,$ there is $k_\a\in\bbbf\setminus\{0\}$ such that $f_\a=k_\a e_\a.$  Similarly, there is a nonzero scalar $k'_\a$ such that $y_\a=k'_\a x_{f(\a)}.$ Define $\phi:\la R\ra\longrightarrow\bbbf\setminus\{0\}$ mapping $\a\in\Pi$ to $k^{-1}_\a$ and $\phi':\la S\ra\longrightarrow\bbbf\setminus\{0\}$ mapping $f(\a)\in f(\Pi)$ to $k'_\a.$ Now  using (\ref{final5}), one can see that isomorphism  $\psi:=\tilde \phi'\circ\theta\circ \tilde \phi$ has  the desired properties. Next suppose  $\Pi$ is a base for $R$ and $\psi$ and $\psi'$ are two isomorphisms from $\gg$ to $\LL$ mapping $f_\a$ to $y_\a$ and $h_\a$ to $h'_{f(\a)}$ for all $\a\in \Pi.$ Then $\varphi:=\psi^{-1}\circ\psi'$ is an automorphism of $\gg$ mapping $f_\a$ to $f_\a$ and $h_\a$ to $h_\a,$ for $\a\in\Pi.$ Since  for $\a\in\Pi,$ $\varphi(f_\a)=f_\a,$ we have $\varphi(e_\a)=e_\a.$ On the other hand, as $\varphi$ is identity on $\hh,$  $\varphi$ preserves the root spaces. This together with the fact that   $[e_\a,\varphi(e_{-\a})]=[\varphi(e_\a),\varphi(e_{-\a})]= \varphi(h_\a)=  h_\a$ for all $\a\in \Pi,$ implies that $\varphi(e_{-\a})=e_{-\a}.$ Now suppose that $\a\in R^\times,$ since $\Pi$ is a base, there are $r_1,\ldots,r_n\in\{\pm1\}$ and $\a_{i_1},\ldots,\a_{i_n}\in \Pi$ such that $\a=r_1\a_{i_1}+\cdots+r_n\a_{i_n}$ and that  $r_1\a_{i_1}+\cdots+r_t\a_{i_t}\in\rcross$ for all $1\leq t\leq n.$ This together with Lemma 2.4 of  \cite{long-ver} and the fact that root spaces corresponding to nonzero roots are one dimensional, implies that $\gg^\a=\bbbf[e_{r_n\a_n},\ldots,[e_{r_3\a_3},[e_{r_2\a_2},e_{r_1\a_1}]]\cdots].$ But $\varphi(e_\a)=e_\a$ and $\varphi(e_{-\a})=e_{-\a}$ for all $\a\in \Pi,$ so the restriction of $\varphi $ to $\gg^\a$ is identity. This completes the proof.
\end{proof}

\subsection{Classification}
\begin{proposition}\label{extention}
Suppose that $(\gg,\fm,\hh)$ and $(\LL,\fm',T)$ are two infinite dimensional  locally finite basic classical simple Lie superalgebras with corresponding root systems $R$ and $S$ respectively. Denote the induced forms on $\v:=\hbox{span}_\bbbf R$ and $\u:=\hbox{span}_\bbbf S$ again by $\fm$ and $\fm'$ respectively. Suppose that $(\la R\ra,\fm,R)$ and $(\la S\ra,\fm',S)$ are isomorphic locally finite root supersystems, then there is an isomorphism from $\gg$ to $\LL$ mapping $\hh$ onto $T.$
\end{proposition}
\begin{proof}
Locally finite basic classical simple Lie superalgebras with zero odd part are exactly locally finite split simple Lie algebras in the sense of \cite{NS}. So contemplating \cite[Thm. IV.4]{NS}, we may assume $\LL_{1}\neq \{0\}.$ Suppose that $f:\la R\ra\longrightarrow \la S\ra$ is the isomorphism from $R$ to $S$ with  $(f(\a),f(\b))'=k(\a,\b)$ for $\a,\b\in R$ and some  nonzero scalar $k.$
Using Lemma \ref{lin2}, we extend $f$ to a linear isomorphism, denoted again by $f,$ from $\hbox{span}_\bbbf R$ to  $\hbox{span}_\bbbf S$ with $(f(\a),f(\b))'=k(\a,\b)$ for $\a,\b\in \hbox{span}_\bbbf R.$
Fix $r,s\in\bbbf\setminus\{0\}$ such that $r=sk.$ Therefore, we have
$$
r(\a,\b)=sk(\a,\b)=s(f(\a),f(\b))'\;\;\;\;\; (\a,\b\in \hbox{span}_\bbbf R).
$$
For $\a\in \hbox{span}_\bbbf R,$ take $t_\a$ to be the unique element of $\hh$ representing $\a$ through $\fm$ and for $\b\in \hbox{span}_\bbbf S,$ take $t'_\b$ to be the unique element of $T$ representing $\b$ through $\fm'.$ Set
$$
h_\a:=rt_\a\andd h'_\b:=st'_\b\;\;\;\;(\a\in R,\;\; \b\in S).
$$
By  Lemma \ref{base}, there is a base $\Pi$ for $R$ and a class $\{R_\gamma\mid \gamma\in\Gamma\}$ of finite irreducible  closed sub-supersystems of $R$ of the same type as $R$ such that $R$ is the direct union of $R_\gamma$'s and for each $\gamma\in\Gamma,$ $\Pi\cap R_\gamma$ is a base for $R_\gamma.$
Now  $\Pi':=f(\Pi)$ is  a base for $S$ and  $S=\cup_{\gamma\in \Gamma} S_\gamma$ in which $S_\gamma:=f(R_\gamma)$ is a finite irreducible  closed sub-supersystem of $S$ and $\Pi'_\gamma:=\Pi'\cap S_\gamma$ is  a base for $S_\gamma.$  For each $\gamma\in\Gamma,$ set $$\gg(\gamma):=\sum_{\a\in R_\gamma^\times}\gg^\a\op\sum_{\a\in R_\gamma^\times}[\gg^\a,\gg^{-\a}]\andd \hh(\gamma):=\sum_{\a\in R_\gamma^\times}[\gg^\a,\gg^{-\a}]$$ and
$$\LL(\gamma):=\sum_{\a\in S_\gamma^\times}\LL^\a\op\sum_{\a\in S_\gamma^\times}[\LL^\a,\LL^{-\a}]\andd T(\gamma):=\sum_{\a\in S_\gamma^\times}[\LL^\a,\LL^{-\a}].$$
Then as in the proof of  Lemma 2.28 of \cite{long-ver}, $\gg(\gamma)$ is  a finite dimensional  basic classical simple Lie superalgebra with Cartan subalgebra $\hh(\gamma)$ and corresponding root system $R_\gamma$ and  $\LL(\gamma)$ is  a finite dimensional  basic classical simple Lie superalgebra with Cartan subalgebra $T(\gamma)$ and corresponding root system $S_\gamma.$ Now fix $\{f_\a\in\gg^\a\mid\a\in \Pi\}$ and $\{y_\a\in\LL^{f(\a)}\mid \a\in \Pi\},$ By Lemma \ref{cor}, for each $\gamma\in\Gamma,$ there is a unique isomorphism $\theta_\gamma$ from $\gg(\gamma)$ to $\LL(\gamma)$ mapping $f_\a$ to $y_\a$ and $h_\a$ to $h'_{f(\a)}$ for $\a\in \Pi_\gamma.$   Now for $\gamma_1,\gamma_2\in\Gamma$ with $R_{\gamma_1}\sub R_{\gamma_2},$ $\theta_{\gamma_1},\theta_{\gamma_2}\mid_{\gg(\gamma_1)}$ are isomorphisms from $\gg(\gamma_1)$ to $\LL(\gamma_1)$ mapping $f_\a$ to $y_\a$ and $h_\a$ to $h'_{f(\a)}$ for all $\a\in \Pi_\gamma,$ therefore, we have $\theta_{\gamma_1}=\theta_{\gamma_2}\mid_{\gg(\gamma_1)}$ by Lemma \ref{cor}.  This allows us to define the isomorphism $\theta:\gg\longrightarrow \LL$ by $\theta(x)=\theta_\gamma(x)$ if $x\in\gg(\gamma).$ The isomorphism  $\theta$ maps $\hh$ onto $T.$ This completes the proof.
\end{proof}

\begin{corollary}
\label{conjugate}
Suppose that $\LL$ is a locally finite basic classical simple Lie superalgebra. Assume  $\hh$ and  $T$ are two Cartan subalgebras of $\LL$ with corresponding  root systems $R$ and $S$ respectively. Then $\hh$ and $T$ are conjugate under $Aut(\LL)$ if and only if  $R$ and $S$ are isomorphic locally finite root supersystems.
\end{corollary}
\begin{proof}
Assume $\phi:\LL\longrightarrow \LL$ is a Lie superalgebra automorphism such that $T=\phi(\hh).$  Define the bilinear form $\fm':\LL\times\LL\longrightarrow \bbbf$ with  $(x,y)':=(\phi^{-1}(x),\phi^{-1}(y))$ for all $x,y\in \LL.$  The linear isomorphism $\phi\mid_{\hh}:\hh\longrightarrow T$ induces a linear isomorphism $\dot\phi:\hh^*\longrightarrow T^*$ mapping $\a\in\hh^*$ to $\a\circ(\phi|_T)^{-1}.$ For each $h\in\hh,$ $\a\in R$ and $x\in\LL^\a, $ we have  $[\phi(h),\phi(x)]=\phi[h,x]=\a(h)\phi(x)=\dot\phi(\a)(\phi(h))\phi(x).$ Now it follows that   $(\LL,\fm',T)$ is an extended affine Lie superalgebra, $\phi(\LL^\a)=\LL^{\dot\phi(\a)}$ and that $S=\dot\phi(R).$  For $\b\in \hbox{span}_\bbbf S,$ take $t'_\b$ to be the unique element of $T$ representing $\b$ through $\fm'$ and for $\a\in \hbox{span}_\bbbf R,$ take $t_\a$ to be the unique element of $\hh$ representing $\a$ through $\fm.$ One can easily check that $\phi(t_\a)=t'_{\dot\phi(\a)}$  and  $(t_\a,t_\b)=(t'_{\dot\phi(\a)},t'_{\dot\phi(\b)})'$ for all $\a,\b\in R.$ These all together imply that $\dot\phi\mid_{\la R\ra}$ defines an isomorphism from $R$ to $S.$ The reverse part follows from Proposition \ref{extention}, \cite[Cor. IV.5]{NS} and finite dimensional theory of Lie superalgebras.
\end{proof}

Using \cite[Thm. IV.6]{NS}, one knows the classification of locally finite split simple Lie algebras, i.e., locally finite basic classical simple Lie superalgebras with zero odd part. In what follows using  Theorem \ref{main1}, Examples \ref{exa1}, \ref{exa2},  Proposition \ref{extention} and  Lemma \ref{class1},  we  give the classification of locally finite basic classical simple Lie superalgebras with nonzero odd part:

\begin{theorem} Each locally finite basic classical simple Lie superalgebra with nonzero odd part is either a finite dimensional basic classical simple Lie superalgebra or   isomorphic to one  and only one of the Lie superalgebras $\mathfrak{osp}(2I,2J)$ ($I,J$   index sets with  $|I\cup J|=\infty,$ $|J|\neq 0$),  $\mathfrak{osp}(2I+1,2J)$ ($I,J$   index sets with  $|I|<\infty,$ $|J|=\infty$) or $\mathfrak{sl}(I_{0},I_{1})$  ($I$ an infinite superset with $I_{0},I_{1}\neq\emptyset$).
\end{theorem}

\begin{proposition}Suppose that $\LL$ is an infinite dimensional  locally finite basic classical simple Lie superalgebra with nonzero odd part, then  if for  an infinite index set  $I$ and a nonempty index set $J,$ $\LL\simeq \mathfrak{osp}(2I+1,2J)\simeq\mathfrak{osp}(2I,2J),$  there are two conjugacy classes for Cartan subalgebras of $\LL$ under  $Aut(\LL);$ otherwise  all Cartan subalgebras of $\LL$ are conjugate under  $Aut(\LL),$ i.e., there is just one conjugacy class  for Cartan subalgebras of $\LL$ under  $Aut(\LL).$
\end{proposition}

\begin{proof}We first assume $I$ is an infinite index set, $J$ a nonempty index set and $\LL\simeq \mathfrak{osp}(2I+1,2J)\simeq\mathfrak{osp}(2I,2J).$ We know form Example \ref{exa1} that there are Cartan subalgebras $\hh_1$ and $\hh_2$ for $\LL$ with corresponding root systems $R_1$ of type $B(I,J)$ and $R_2$ of type $D(I,J)$  respectively; in particular thanks to Corollary \ref{conjugate}, there are at least two conjugacy classes for Cartan subalgebras of $\LL$ under $Aut(\LL).$ We next note that there is    a decomposition  $\LL_{0}=\gg^1\op\gg^2$ for $\LL_{0}$ into simple ideals in which  $\gg^1$ is isomorphic to $\mathfrak{o}(2I+1,\bbbf)\simeq \mathfrak{o}(2I,\bbbf)$ and  $\gg^2$  is isomorphic to $\mathfrak{sp}(J,\bbbf);$ see \cite{NS} for the notations.  By \cite[Cor. VI.8]{NS} and finite dimensional theory of Lie algebras, there are two  conjugacy classes for Cartan subalgebras of $\gg^1$ under $Aut(\gg^1)$ and there is just one conjugacy class for Cartan subalgebras of $\gg^2$ under $Aut(\gg^2).$
Therefore, up to $Aut(\gg^1)$-conjugacy,  $\hh_1\cap \gg^1,\hh_2\cap \gg^1$ are the only non-conjugate  Cartan subalgebras of $\gg^1;$ also  $\hh_1\cap \gg^2,\hh_2\cap \gg^2$  are $Aut(\gg^2)$-conjugate Cartan subalgebras of $\gg^2$  and in fact  up to $Aut(\gg^2)$-conjugacy,  $\hh_1\cap \gg^2$ is  the only  Cartan subalgebra of $\gg^2.$  Now suppose that $T$ is a Cartan subalgebra of $\LL$ with corresponding root system $S=S_0\cup S_1.$ We want  to show that $T$ is either conjugate to $\hh_1$ or to $\hh_2.$ Since $T\cap \gg^1$ is a Cartan subalgebra of $\gg^1$ and $T\cap \gg^2$ is a Cartan subalgebra of $\gg^2,$ there are $i\in\{1,2\}$ and  $\phi_1\in Aut(\gg^1),\phi_2\in Aut(\gg^2)$ such that $\phi_1(T\cap \gg^1)=\hh_i\cap\gg^1$ and $\phi_2(T\cap \gg^2)=\hh_i\cap\gg^2.$ So $\phi_1\op\phi_2$ is an automorphism of $\LL_{0}$ mapping $T=(T\cap \gg^1)\op(T\cap \gg^2)$ to $\hh_i=(\hh_i\cap \gg^1)\op(\hh_i\cap \gg^2).$
 This implies that $(R_i)_0$ is isomorphic to $S_0.$ So using the  classification  of locally finite root supersystems (Theorem \ref{classification}) together with  Proposition 2.5 and  Lemma 2.17 of \cite{long-ver} and the fact that $|R_i|,|S|=\infty,$ $R_i$ is isomorphic to $S.$
Therefore there is an automorphism of $\LL$ mapping $T$ to $\hh_i$ by proposition \ref{extention}.  This implies that there are exactly two conjugacy classes for Cartan subalgebras of $\LL$ under $Aut(\LL).$

Next suppose that  $\LL $ is one of the Lie superalgebras  $ \mathfrak{osp}(2I,2J),\mathfrak{osp}(2I+1,2J)$ where $I,J$ are   index sets with  $0\neq|I|<\infty,$ $|J|\neq 0$ or   $\mathfrak{sl}(I_{0},I_{1})$ where   $I$ is an infinite superset with $|I|=\infty$ and  $I_{0},I_{1}\neq\emptyset.$
Take $\fh$ to be the standard Cartan subalgebra of $\LL$ introduced in  Examples \ref{exa1} and \ref{exa2} and consider its corresponding root system $R.$ Next suppose that $T$ is another Cartan subalgebra of $\LL$ and take  $S$ to be the corresponding root system of $\LL$ with respect to $T.$ Then $S$ is an irreducible locally finite root supersystem with $|S|=\infty.$ From Theorem \ref{classification} and Lemmas 2.16 and 2.17, $S$ is isomorphic to the root system of one of  the Lie superalgebras $\mathfrak{a}_{I',J'},\mathfrak{b}_{I',J'},\mathfrak{c}_{T'}$ introduced in Lemma \ref{class1}. Call this Lie superalgebra $\gg$ and take $\hh$ to be its standard Cartan subalgebra, so by Proposition \ref{extention}, there is an isomorphism  $\phi:\LL\longrightarrow  \gg$ mapping $T$ to $\hh.$  Now since $\LL\simeq\gg,$ using Lemma \ref{class1}, there is an isomorphism $\psi$ form $\LL$ to $\gg$ mapping $\fh$ to $\hh.$ Therefore, $\psi^{-1}\circ \phi$ is an automorphism of $\LL$ mapping $T$ to $\fh.$ This completes the proof.\end{proof}

\centerline{\bf Acknowledgment}
This research was in part
supported by a grant from IPM (No. 93170415) and partially carried out in IPM-Isfahan branch. The author acknowledges this support.

\end{document}